\documentclass[a4paper, 10pt]{article}
\usepackage{yfonts}
\usepackage{comment}
\usepackage[normalem]{ulem}
\usepackage{fancyhdr}
\usepackage{bm}
\usepackage{cancel}
\usepackage{latexsym}
\usepackage{color}
\usepackage{cancel}
\usepackage{amssymb,amsmath,amsfonts,bbm,pifont,upgreek,bbold
}

\usepackage[nointlimits]{esint}
\usepackage[colorlinks=true]{hyperref}
\hypersetup{urlcolor=blue, citecolor=red}
\usepackage{graphicx}
\usepackage{tikz}
\newcommand\beq[1]{\begin{equation}\label{#1} }

\renewcommand{\theequation}{\arabic{section}.\arabic{equation}}

\newtheorem{theorem}{Theorem}[section]
\newtheorem{definition}{Definition}[section]
\newtheorem{proposition}{Proposition}[section]
\newtheorem{lemma}{Lemma}[section]

\newtheorem{sublemma}{Sublemma}[section]
\newtheorem{remark}{Remark}[section]
\newtheorem{notationalremark}{Notational Remark}[section]
\newtheorem{corollary}{Corollary}[section]
\newtheorem{assumption}{Assumption}[section]
\newtheorem{claim}{Claim}[section]

\newtheorem{problem}{Problem}[section]
\newtheorem{tools}{$\negsp\negsp$}[subsection]

\newcommand\thm[1]{\begin{theorem}\label{#1}}
\newcommand\thmtwo[2]{\begin{theorem}[#1]\label{#2}}
\newcommand\ethm{\end{theorem} }
\newcommand\dfn[1]{\begin{definition}\label{#1} \rm}
\newcommand\dfntwo[2]{\begin{definition}[#1]\label{#2} \rm}
\newcommand\edfn{\end{definition} }
\newcommand\pro[1]{\begin{proposition}\label{#1}}
\newcommand\protwo[2]{\begin{proposition}[#1]\label{#2}}
\newcommand\epro{\end{proposition} }
\newcommand\lem[1]{\begin{lemma}\label{#1}}
\newcommand\lemtwo[2]{\begin{lemma}[#1]\label{#2}}
\newcommand\elem{\end{lemma} }
\newcommand\sublem[1]{\begin{sublemma}\label{#1}}
\newcommand\sublemtwo[2]{\begin{sublemma}[#1]\label{#2}}
\newcommand\esublem{\end{sublemma} }
\newcommand\rem[1]{\begin{remark}\label{#1} \rm}
\newcommand\erem{\end{remark} }
\newcommand\notrem[1]{\begin{notationalremark}\label{#1} \rm}
\newcommand\enotrem{\end{notationalremark} }
\newcommand\cor[1]{\begin{corollary}\label{#1}}
\newcommand\cortwo[2]{\begin{corollary}[#1]\label{#2}}
\newcommand\ecor{\end{corollary} }
\newcommand\asmp[1]{\begin{assumption}\label{#1}}
\newcommand\asmptwo[2]{\begin{assumption}[#1]\label{#2}}
\newcommand\easmp{\end{assumption} }
\newcommand\clm[1]{\begin{claim}\label{#1}}
\newcommand\eclm{\end{claim} }

\expandafter\chardef\csname pre amssym.def
at\endcsname=\the\catcode`\@
\catcode`\@=11
\def\undefine#1{\let#1\undefined}
\def\newsymbol#1#2#3#4#5{\let\next@\relax
 \ifnum#2=\@ne\let\next@\msafam@\else
 \ifnum#2=\tw@\let\next@\msbfam@\fi\fi
 \mathchardef#1="#3\next@#4#5}
\def\mathhexbox@#1#2#3{\relax
 \ifmmode\mathpalette{}{\m@th\mathchar"#1#2#3}%
 \else\leavevmode\hbox{$\m@th\mathchar"#1#2#3$}\fi}
\def\hexnumber@#1{\ifcase#1 0\or 1\or 2\or 3\or 4\or 5\or 6\or 7\or
8\or
 9\or A\or B\or C\or D\or E\or F\fi}
\ifcase\@ptsize
 \font\tenmsb=msbm10
 \font\sevenmsb=msbm7
 \font\fivemsb=msbm5
\or
 \font\tenmsb=msbm10 scaled \magstephalf
 \font\sevenmsb=msbm7 scaled \magstephalf
 \font\fivemsb=msbm5  scaled \magstephalf
\or
 \font\tenmsb=msbm10 scaled \magstep1
 \font\sevenmsb=msbm7 scaled \magstep1
 \font\fivemsb=msbm5 scaled \magstep1
\fi
\newfam\msbfam
\textfont\msbfam=\tenmsb
\scriptfont\msbfam=\sevenmsb
\scriptscriptfont\msbfam=\fivemsb
\edef\msbfam@{\hexnumber@\msbfam}
\def\Bbb#1{\fam\msbfam\relax#1}
\def\widehat#1{\setboxz@h{$\m@th#1$}%
 \ifdim\wdz@>\tw@ em\mathaccent"0\msbfam@5B{#1}%
 \else\mathaccent"0362{#1}\fi}
\def\widetilde#1{\setboxz@h{$\m@th#1$}%
 \ifdim\wdz@>\tw@ em\mathaccent"0\msbfam@5D{#1}%
 \else\mathaccent"0365{#1}\fi}

\def\RIfM@{\relax\ifmmode}
\def\nonmatherr@#1{\errmessage{\string#1\space allowed only in math mode}}
\def\Bbb{\RIfM@\expandafter\Bbb@\else
 \expandafter\nonmatherr@\expandafter\Bbb\fi}
\def\Bbb@#1{{\Bbb@@{#1}}}
\def\Bbb@@#1{\fam\msbfam\relax#1}
\def\setboxz@h{\setbox\z@\hbox}
\def\wdz@{\wd\z@}
\catcode`\@=\csname pre amssym.def at\endcsname
\font\tengrm=gsmn1000
\font\ninegrm=gsmn0900

\font\sevengrm=gsmn0700

\font\fivegrm=gsmn0500
\newfam\grmfam
 \def\grm{\fam\grmfam\tengrm} \textfont\grmfam=\tengrm
\scriptfont\grmfam=\sevengrm \scriptscriptfont\grmfam=\fivegrm
\def\grm{\fam\grmfam\ninegrm} \textfont\grmfam=\ninegrm
 \scriptfont\grmfam=\sevengrm \scriptscriptfont\grmfam=\fivegrm

\newcommand{\ii}{{\rm i}  }

\newcommand{\negsp}{\hspace{-.09truecm}}

\newcommand\sign{{\, \rm sign\, }}

\renewcommand{\natural}{{\Bbb N}   }

\newcommand{{\cE}}{{\cal  E} }

\newcommand{{\cH}}{{\cal H} }
\newcommand{{\cK}}{{\cal K} }

\newcommand{{\cJ}}{{\cal J}}

\usepackage{color}
\definecolor{yellow}{rgb}{0.99, 0.93, 0.0}

\def\norm#1{\left |#1\right |}
\def\Norm#1{\left \|#1\right \|}

\begin{document}
\title{
Improved stability estimates at elliptic  equilibria of Hamiltonian systems\vskip.1in\large
}

\author{Massimiliano Guzzo, Chiara Caracciolo, Gabriella Pinzari
\\ \footnotesize  Universit\`a degli Studi di Padova 
\\ \footnotesize  Dipartimento di Matematica ``Tullio Levi-Civita''
\\ \footnotesize  Via Trieste 63 - 35121 Padova, Italy 
}

\maketitle

\begin{abstract}\footnotesize{
{This paper deals with an improvement of the ``a--priori stability bounds''
   on the variation of the action variables and on the stability time
   obtained from a  given Birkhoff normal form around the elliptic equilibrium  point of an Hamiltonian  system satisfying a non-resonance condition of
   finite order $N$. In particular,  we improve the standard a--priori lower bound on the stability time from a purely linear dependence on the
   inverse of the polynomial norm of the remainder of the normal form to the sum of a linear  term
   (which is still present but with a different constant coefficient) and a quadratic one. 
   The prevalence between the linear and the quadratic term depends on the resonance
   properties of all the monomials in the remainder of the normal form with degree from $N$ to a finite order $M$. We also provide a comparative example of the new estimates and the traditional a priori ones in the framework of computer-assisted proofs.
 } 
}
\end{abstract}

\maketitle

\tableofcontents

\renewcommand{\theequation}{\arabic{equation}}
\setcounter{equation}{0}

\section{Introduction}\label{Purpose of the paper}
        This paper deals with an improvement on the ``a--priori lower bounds''
        on the variation of the action-variables and on the stability time
 around the elliptic equilibrium  of an Hamiltonian  system 
   obtained from a  given Birkhoff normal form. The setting is as follows.
  \vskip.1in
 \noindent
 We consider a Hamilton function  
\begin{eqnarray}
  H_0(p, q)=\sum_{k\ge 2}h_k(p, q)
  \label{H0}
\end{eqnarray}
having the form of a convergent power series in a neighborhood $D_0$ of $(0,0)\in {\mathbb R}^n\times {\mathbb R}^n$ 
where  $(p,q)$
are canonical variables\footnote{The $q=(q_1\,,\ldots\,, q_n)$ are ``coordinates'' and the $p=(p_1\,,\ldots\,, p_n)$ are conjugate momenta, namely 
they evolve accordingly to 
$$\displaystyle\left\{
\begin{array}{ll}
\dot q_i=\{q_i,H_{0}\}\\
\dot p_i=\{p_i,H_{0}\}
\end{array}
\right.$$
with $\{f,g\}=\nabla_q f \cdot \nabla_p g-\nabla_q g \cdot \nabla_p f$.} and the functions $h_k$  are homogeneous polynomials of degree $k$, with $h_2(p, q)$ having the form
\begin{eqnarray}
\displaystyle h_2(p, q)=\sum_{j=1}^n\Omega_{0j} \frac{p_j^2+q_j^2}{2}\,.
\end{eqnarray}
When the components of the frequency vector $\Omega_0=(\Omega_{01}, \ldots, \Omega_{0n})$ have different signs the application of  standard Lyapunov stability theorem (which would grant perpetual
stability around the origin) is prevented, and so the determination of an a--priori lower bound of the  stability time is usually posed as follows.
 \begin{problem}\label{problem of elliptic equilibrium} For any 
   $D_1\subset D_0$ neighborhood of   $(p,q)=(0,0)$,
   provide a neighborhood $D_2\subseteq D_1$ of  $(p,q)=(0, 0)$ and 
   a lower bound of the time $T$ such that any solution $(p(t), q(t))$  of the Hamilton equation of $H_0$ with initial datum $(p(0), q(0))\in D_2$
   stays in $D_1$ for all $t$, with $|t|\le T$.
 \end{problem}

 If the  vector of the frequencies $\Omega_{0}$
 satisfies the non-resonance condition of finite order $N$:
 \begin{eqnarray}
  \Omega_0\cdot k  \ne 0\quad \forall\ k\in {\mathbb Z}^n\,: \quad 0<\sum_{j=1}^n|k_j|\le N\,,
  \label{non-resN}
\end{eqnarray}
Hamiltonian  (\ref{H0})  is conjugate by a close to the identity canonical
 transformation to the so-called Birkhoff normal form \cite{birkhoff1927}: 
 \begin{equation}
   {\cal H}_N = H_N(p,q)+f_N(p,q)
   \label{hfN}
\end{equation}   
where $f_N(p,q)$ is a convergent series of monomials in the variables $q,p$
of degree larger than $N$, and the function $H_N$ 
depends on the variables $q,p$ only through the  ``action functions''
\begin{eqnarray}\label{action functions}
  I_j(p, q)=\frac{p_j^2+q_j^2}{2}\,,
\end{eqnarray}
and so in particular we have $\{ I_j , H_N \} = 0$. In the following, we also refer to the polynomial function
\begin{equation}\label{Zfunction}
  Z(I)= Z_2(I)+\ldots+ Z_{2s}(I)
\end{equation}
where $Z_{2j}(I)$ are homogeneous components with degree $j$  in the indeterminates $I=(I_1, \ldots, I_n)$, such that:
\begin{equation}
Z(I(p, q))=H_N(p,q)  .
\end{equation}
When $n=2$, and the Birkhoff normal form of some order $N$ satisfies
a non--degeneracy condition, the perpetual stability
of the equilibrium point may be recovered using KAM theory \cite{arnold64,Arnold-eq,Russmann:1976,meyer-schimdt}). 
Instead, when $n\geq 3$, although KAM invariant tori provide perpetual stability in a large measure subset of a neighbourhood of the origin, they 
may coexist with instability since their distribution may not prevent the
possibility for the action variables to afford arbitrary variations in finite times. In these cases, the stability of the equilibrium point and of the action variables have been studied over finite time intervals $\norm{t}<T$. 

\begin{remark}\rm
Research on KAM theory has significantly advanced our understanding of dynamics near equilibrium points, particularly in addressing the crucial question of long-term stability in planetary systems ($n+1$ bodies mutually
interacting through gravity only). The proof of KAM stability around the so-called ``secular equilibrium'' configuration -- consisting of $n$ planets describing co-planar and co-circular orbits about the Sun -- has been obtained in a series of papers \cite{arnold63,laskarR95,fejoz04,locatelligiorgilli2000,chierchiaPi11b,FH2025}.
On the side of instability, lower-dimensional hyperbolic KAM tori, coexisting
with the fully dimensional ones,  have been found in \cite{pinzariL2023} for the
three--body problem and Arnold diffusion  has been proved in
\cite{clarkeFG2024}  for a four--body case. Remarkable
numerical examples of instability of the secular equilibrium configuration
of our Solar System 
have been provided in \cite{laskar94,laskar08,mogavero-hoang-laskar}. 
\end{remark}

Along the vein of Nekhoroshev theory \cite{nekhoroshev77,nekhoroshev79,benettinG1986,lochak1992,poschel93,lochakNN1994,bonemouraN12,guzzoCB16}, exponentially-long  times  of  stability of the elliptic equilibrium have been obtained 
under additional assumptions on $Z(I)$, see Remark \ref{remark2} \cite{fassoGB1998,guzzoFB1998,BenettinFassoGuzzo1998,niederman1998,niederman2013,bounemouraFN2020nonl,bounemouraFN2020,palacianetal2021}. A common strategy to obtain a lower bound for the stability time $T$ in Problem \ref{problem of elliptic equilibrium} is to resort to  a ``confinement of actions'' for all $\norm{t}<T$ and for all the initial conditions in the sets:
\begin{equation}\label{setD2}
B(R)=\{(p,q): |I_1(p,q)|,\ldots ,|I_n(p,q)|< R^2\}
\end{equation}
which is valid for any arbitrary small $R$.  Indeed, for suitably chosen $c,R_0> 0$,
the objective is to prove that for all  $R<R_0$, and all the solutions  of the Hamilton equations of (\ref{hfN}) with initial conditions in $B(R)$, one has: 
\begin{eqnarray}\label{confinement of actions}
  && |I_j(t)-I_j(0)|<cR^2\,,\quad \forall\ j=1\,,\ldots\,,n
\quad \forall \ t:\ |t|<T  .
\end{eqnarray}
The time  $T$  in \eqref{confinement of actions} is usually referred to as ``actions confinement time'' and,
when the constant $c$ is suitably small, the confinement on
the actions provides a stronger stability result than the one stated in Problem 1
(see Section \ref{thestatements} for details and precise statements).
\vskip 0.4 cm

The question we consider in this paper is to provide the longest possible
quantitative lower bound for $T$ when no further normalization at orders larger than $N$ is available, and no other hypotheses on the Hamiltonian (\ref{H0}) are
assumed (in particular, we do not use any non-degeneracy or steepness
condition on $Z(I)$).  
\begin{remark}\rm
We may have this situation in two cases:  
\begin{itemize}
\item[-]  Even if the frequency vector $\Omega_0$
  is non-resonant (for example satisfying a Diophantine condition)  one cannot expect to improve the stability bounds by defining a normal form
  of arbitrarily large order $N$. 
Indeed, in general, the remainder term $f_N$ in \eqref{hfN} does not go to zero in a given set $B(R)$ as $N\to \infty$. Rather,
the polynomial norm
of $f_N$ in any given set $B(R)$ may be estimated  from above by a sequence
which monotonically decreases with $N$ only up to an optimal order $N_{opt}$ \cite{Moser1,Littlewood1,Littlewood2,giorgillietal1989}.  Therefore, Birkhoff normal forms of order $N$ larger than $N_{opt}$ do not provide any improvement on the stability estimates. 
  
\item[-] The Birkhoff normal forms of a given Hamiltonian
  (\ref{H0}) can be computed by iterative algorithms using the rigorous approach of the interval arithmetics (see \cite{KSW96,Tucker15} for an introduction to validated numerics and applications to analysis): on the one hand, all the coefficients of the monomials of the normal forms  up to some degree $M$ (typically $M > N$) are provided within intervals; on the other hand the tail for the remainder of the Taylor series of degree larger than $M+1$ is rigorously estimated
  using standard analytic methods (see, e.g., \cite{giorgillietal1989,celgio1991,carloc2020},  where applications to the restricted three-body problem are considered). Therefore, the Birkhoff normal forms which are computed this way allow to provide rigorous stability estimates for various problems which are important for applications, and for which the purely analytic techniques do not provide sufficiently long stability times \cite{giorgilliLS09,sansotteraLG13,GLS2017}. For these problems, the goal is to prove the existence of normal forms of suitably large order $N$ that, although it may not be the optimal order, allows to prove a stability time which is sufficiently long for the specific problem (for example, in the case of a planetary problem, the age of the universe would be
  sufficiently long).
  \end{itemize}
\end{remark}

\begin{remark}\label{remark2}
  \rm If further hypotheses are assumed, such as quasi-convexity, steepness
  or directionally quasi-convexity of the Birkhoff normal form,
  specific Nekhoroshev's stability results proved for elliptic equilibria apply \cite{fassoGB1998,guzzoFB1998,BenettinFassoGuzzo1998,niederman1998,niederman2013,bounemouraFN2020nonl,bounemouraFN2020,palacianetal2021}. Nekhoroshev's estimates
  provide stability times which increase exponentially, or super-exponentially,
  as the distance of the initial condition from the equilibrium
  goes to zero, and are proved using specific resonant or non-resonant
  normal forms. The argument that we specifically develop in this paper for the Birkhoff normal form can be used also for any type 
  of non-resonant normal form defined within the Nekhoroshev's theory
  (more details will be given in Section \ref{argumentforAA}),  while its rigorous extension to the resonant normal forms is still open
  (a semi-analytic treatment of resonances of multiplicity 1 is   provided in \cite{guzzoEP2020}).
\end{remark}

The inequality \eqref{confinement of actions} is typically managed
in the literature of Nekhoroshev theory as an a--priori bound provided directly by Hamilton equations of Hamiltonian \eqref{hfN}, representing a non--resonant normal form of the original Hamiltonian system (\ref{H0}): for any orbit which does not leave the ball $B(R)$ along the interval $[0,t]$,  one has
    \begin{eqnarray}\label{a priori boundOLD}
    |I_j(t)-I_j(0)|=\left|\int_0^t\{I_j, f\}d\tau\right|\le
  C R^{N+1} |t|
    \end{eqnarray}
    where $C$ is determined by $f_N$. The inequality
    \eqref{a priori boundOLD} implies \eqref{confinement of actions}
    for all initial conditions in $B(R/2)$ and 
    for all $\norm{t}\leq T$ with
 \begin{eqnarray}\label{a priori bound}T\ge c_* R^{-({N}-1)}\end{eqnarray}
 Inequality \eqref{a priori bound} represents the  state--of--the--art  for what concerns the lower bounds for the stability time $T$
 for the case $n\geq 3$, 
if no other assumption is taken up on $\Omega_0$ (e.g., we are not assuming $\Omega_0$ diophantine) or on the Hamiltonian (\ref{H0}) (see, for example, Remark
\ref{remark2}). 
\vskip 0.4 cm

In this paper we provide an improvement of inequality (\ref{a priori boundOLD}) 
by developing in a rigorous way a recent point of view taken up in
\cite{guzzoEP2020},  where {semi--analytic} bounds on the speed of diffusion have been obtained using arguments of stationary and non-stationary phase analysis of the remainder terms. The semi-analytic argument presented in \cite{guzzoEP2020} is here considered and made rigorous in the context of
the finitely non-resonant elliptic equilibria. 
To explain our main argument, let us first introduce the complex Birkhoff canonical variables (see \cite{birkhoff1927}):
\begin{eqnarray}\label{wz}
w_j=\frac{p_j+\ii q_j}{\ii \sqrt2}\,,\qquad z_j=\frac{p_j-\ii q_j}{\sqrt2}
\end{eqnarray}
and consequently rewrite the Birkhoff normal form ${\cal H}_N$ in Eq. (\ref{hfN}) as
\begin{eqnarray}\label{Hamiltonian}
{\cal H}(z,w)=Z(\ii wz)+ f(z,w)
\end{eqnarray}
with $\ii wz$ denoting the $n$--ple $(\ii w_1z_1\,,\ldots, \ii w_nz_n)$,
and $f$ is a analytic function in the complex ball
  \begin{eqnarray}{{D(R_0)}}:=\{(z,w)  \in 
  {\mathbb C}^n\times{\mathbb C}^n:\ \|(z,w)\|\le R_0\}\end{eqnarray} 
  where
    \begin{eqnarray} \|(z,w)\|:=\max_{i, j}\{|z_i|, |w_j|\}\,,\end{eqnarray} 
    whose Taylor expansion\footnote{For all integer $p$-uples of indices 
    $\nu = (\nu_1,\ldots ,\nu_p)\in {\mathbb Z}^p$   we denote 
    $|\nu|:=\sum_{i=1}^p|\nu_i|$. Notice that we use the same notation
    $\norm{x}$ for the hermitian norm of any $x\in {\Bbb C}^p$.
For any $x\in {\Bbb C}^p$ and $h\in {\Bbb N}^p$ we also
    use the multi-index notation: $x^h=
    x^{h_1}\cdots x^{h_n}$.     } 
   \begin{eqnarray}\label{f}f(z,w)=\sum_{(h, k)\in{\mathbb N}^n\times {\mathbb N}^n\atop{{|h|+|k|}\ge {N{+1}}}}c_{hk} w^h z^k\end{eqnarray}
 starts with degree $N+1$.     
 We limit to consider solutions of the Hamilton equations of $H$
whose initial datum  verify\footnote{The ``bar'' denote complex conjugacy.}
 \begin{eqnarray}\label{reality at time zero}
   \overline w=\ii z\end{eqnarray}
  which correspond to real initial values of the coordinates $(p, q)$ (see \eqref{wz}). As the original Hamiltonian \eqref{H0} is real--valued,
  condition \eqref{reality at time zero} for the initial datum implies
  that it is satisfied along the solution. Therefore, the stability around the equilibrium
 point for the all the solutions $(w(\tau),z(\tau))$ which
belong for all $\norm{\tau}\leq \norm{t}$ to a set:
 \begin{equation}
   D_1 = \{ (z,w): \overline w=\ii z, \, {\rm and}\, \Norm{(z,w)}\leq R\}
 \end{equation}  
 is traditionally obtained from the inequality:
 \begin{eqnarray}\label{a priori boundEq}
    |I_j(t)-I_j(0)| &\le& 
    \sum_{(h,k)\in {\Bbb N}^n\times {\Bbb N}^n}\norm{h_j-k_j}\norm{c_{hk}}
    \left|\int_0^t w(\tau)^hz(\tau)^kd\tau \right| \cr
    &\leq & |t| \, \sum_ {(h,k)\in {\Bbb N}^n\times {\Bbb N}^n}
\norm {c_{hk}} \norm{h_j-k_j}R^{\norm{h}+\norm{k}}   .
 \end{eqnarray}

\vskip 0.4 cm

 We note that in the right--hand side of inequality (\ref{a priori boundEq})
 only monomials of degree $\norm{h}+\norm{k}>N$ are included,
 since the monomials of lower degree appearing in the Birkhoff normal
 form are characterized by $k_j=h_j$ for all $j=1,\ldots ,n$. Given any
 $M>N$, we obtain a better estimate by providing an upper bound for each integral:
\begin{equation}\label{integralwz}
  \left|\int_0^t w(\tau)^hz(\tau)^kd\tau \right|
\end{equation}
as follows:
\begin{itemize}
\item[(i)] if $\norm{h}+\norm{k}>M$ we  use the traditional
  estimate:
  \begin{equation}\label{integralwzI}
  \left|\int_0^t w(\tau)^hz(\tau)^kd\tau \right| \leq R^{\norm{h}+\norm{k}}\norm{t}
\end{equation}
\item[(ii)] if $N<\norm{h}+\norm{k}\leq M$ and
  \begin{equation}
\norm{\sum_{j=1}^n(k_j-h_j)\Omega_{0j}}< a:=4L MR^2
  \end{equation}
  where  $L$ is a   Lipschitz constant for $Z''(I)$ (defined  in (\ref{L})),
  we  use the traditional
  estimate (\ref{integralwzI}) as well. 
\item[(iii)] if $N<\norm{h}+\norm{k}\leq M$ and
  \begin{equation}
\norm{\sum_{j=1}^n(k_j-h_j)\Omega_{0j}}\geq a
  \end{equation}
  we estimate the integral (\ref{integralwz}) by exploiting the
  identity (see Section \ref{The proofs}, Lemma \ref{lemmafirstint}, for a more
  precise formulation and proof):
    \begin{equation}\label{integralwzII}
  \int_0^t w(\tau)^hz(\tau)^kd\tau  
    = {w(\tau)^hz(\tau)^k\over  \ii (h-k) \cdot \Omega (I(\tau)) }
  \Big \vert_0^t 
  - \int_0^t \left \{ {w^h z^k \over  \ii (h-k) \cdot \Omega },f \right \}
  (w(\tau),z(\tau))d\tau ,
\end{equation}
 where:
\begin{eqnarray}\label{Omega0}
  \Omega_j(I):=\partial_{I_j}Z(I)  .
\end{eqnarray} 
\end{itemize}
For the terms of the series satisfying conditions (i) or (ii) we have no improvements with respect to the traditional estimate. Instead, for
the terms satisfying condition (iii),  from Eq. (\ref{integralwzII}) we get indeed a big improvement of inequality (\ref{a priori boundEq}) when $\norm{t}$
  is suitably large. In fact,
both integrals in the right--hand side of Eq. (\ref{integralwzII}) have
integrand of magnitude $R^{\norm{h}+\norm{k}+N-1} \leq R^{2N}$ (i.e. they
are quadratic in $R^N$), while the first term, which is of magnitude $R^{\norm{h}+\norm{k}} \leq R^{N+1}$, is just a boundary term (it is not multiplied by $\norm{t}$). Precise statements about estimates on the variation of the
actions and on the stability time will be provided  Section
\ref{thestatements}.

\begin{remark}\label{remarkM=2N}
{\rm The prevalence between the estimates from the terms satisfying condition
(i), (ii) or (iii) depends on the resonance properties of all the monomials in the remainder of the normal form with degree from $N$ to $M$. For example,
by choosing $M=2N$, the contribution of all
the terms satisfying (i) is dominated by $R^{2N}$, and therefore
is smaller, or comparable, to the contributions from 
the terms satisfying (iii). Therefore, the improvement is very effective when
the terms satisfying (ii) represent a small fraction of the
terms satisfying (iii), a condition which is satisfied as soon as $a$ is small, i.e. for:
\begin{equation}
  R^2 << {1\over 8LN} .
\end{equation}
{For example, for $n=3$, while the number $\#\Lambda_N$ of integer vectors
  $\nu$ with $\norm{\nu}=N$ increases at most with $N^2$, the number $\#\Lambda_{a,N}(\Omega)$ of these vectors with $\norm{\nu \cdot \Omega}<a$ (for $\Omega\ne 0$
  and $a$ suitably small) is in the range between $0$ and a quantity that 
increases at most linearly with $N$. A more
  precise representation of $\#\Lambda_{a,N}(\Omega)$ will be discussed in
Section \ref{An upper bound on the resonant term}.} For the optimal implementation of the argument within a computer assisted proof it is more convenient to keep $M$ as a free parameter, and directly check the optimal choice for the estimate of the
stability time.
}
\end{remark}

\begin{remark}\label{remark-aa}
  {\rm The argument here presented for the Birkhoff normal
    forms can be formulated and made rigorous also for the non--resonant
    normal forms which are defined within Nekhoroshev's theory for
    ``close-to-be integrable, action-angle'' analytic Hamiltonians having
    the form
\begin{equation}
    H_0(I, \varphi)=h(I)+\epsilon f(I, \varphi)
\end{equation}
where $I\in V\subseteq {\mathbb R}^n$, $V$ open and connected, $\varphi\in {\mathbb T}^n$, with ${\mathbb T}:={\mathbb R}/2\pi{\mathbb Z}$;  $\epsilon$ is a very small positive number, see Section \ref{argumentforAA}.
}
\end{remark}

\vskip 0.4 cm

The paper is organized as follows. In Section \ref{thestatements}
we provide the rigorous statements about the bounds on the variation
of action functions and of the stability time around the equilibrium of
a Birkhoff normal form, namely Theorems \ref{main simplified}, \ref{main}
and \ref{main middle}; Sections 3 is devoted to the
proof of the Theorems; in Section 4 we compute
with rigorous arithmetics the estimates provided by Theorems \ref{main simplified}, \ref{main} and \ref{main middle} for a model system, and compare them
with the traditional ones; in Section 4 we discuss the
number of stationary terms that we may have in the remainder terms
of degree between $N+1$ and $M$;  in Section \ref{argumentforAA} we discuss the application of the non-stationary argument to any non-resonant Nekhoroshev normal
form defined with action--angle variables.

\section{The statements}\label{thestatements}

We start with the Hamiltonian (\ref{Hamiltonian}), in Birkhoff
normal form of order $N$:
\begin{eqnarray}
H(z,w)=Z(\ii wz)+ f(z,w)  ,
\end{eqnarray}
and analytic in the complex ball $D(R_0)$. We denote by
\begin{eqnarray}\label{L}
L:=\max_{1\le i\le n}\max_{(z,w)\in {{D(R_0)}}} \sum_{j=1}^n |Z''(\ii wz)_{ij}|
\end{eqnarray}
where $Z''(I)_{ij}:=\partial_{I_j}\Omega_i(I)=\partial^2_{I_iI_j}Z(I)$  is the $n\times n$ Hessian matrix of $Z$. Moreover, we let, for any $0<R\le R_0$,
\begin{eqnarray}\label{Pnorm}\|f\|_{R}:=\sum_{h, k
 }|c_{hk}|R^{|h|+|k|}\end{eqnarray}
 where $c_{hk}$ are the coefficients of the Taylor expansion \eqref{f} of $f(z,w)$.

Finally, if, for any set ${\cal I}\subseteq{\mathbb N}^n\times {\mathbb N}^n$,  $\Pi_{\cal I} f:=\sum_{(h, k)\in{\cal I}  }c_{hk} w^h z^k$ denotes the projection of $f$ on ${\cal I}$, we let
\begin{eqnarray}\label{fff}
f_{0}:=\Pi_{{\cal I}_0}  f\,,\qquad f^*:=\Pi_{{\cal I}^*}  f\,,\qquad f_{*}:=\Pi_{{\cal I}_*}  f
\end{eqnarray}
where, for fixed $a>0$, $M\ge N+1\in \mathbb N$, ${\cal I}_{ 0}$ (respectively, ${\cal I}_{*}$) is the set of $(a, N, M)$--non resonant (respectively, resonant) indices, while 
  ${\cal I}^*$ is the set of $M$--ultraviolet indices:
  \begin{eqnarray}\label{indices}
&& {\cal I}_{ 0}:=\Big\{(h, k):\ N{+1}\le {|h|+|k|}\le M\,,\  |\Omega_0\cdot(h-k)|\ge{a}\Big\}\,,\nonumber\\
&&{\cal I}_{*}:=\Big\{(h, k):\ N{+1}\le {|h|+|k|}\le M\,,\  |\Omega_0\cdot(h-k)|<{a}\Big\}\,,\nonumber\\
&&{\cal I}^*:= \Big\{(h, k):\ {|h|+|k|}> M\Big\}\,.
  \end{eqnarray}
	\vskip.1in
	\noindent
\begin{theorem}\label{main simplified}
  Let $H$, as in \eqref{Hamiltonian}, be in Birkhoff normal form of order $N$ 
and analytic in $D(R_0)$. Let $M\in \mathbb N$, with    $M\ge N{+1}$, 
  $\alpha\in(0, {1/2}]$ and $a>0$. Then, for any  $0<R\le R_0/2$ verifying
\begin{eqnarray}\label{smallness condition}&&\frac{2LM(1+2\alpha)^2 R^2}{{a}}\le 1\\
\label{smallness condition NEW}&&\frac{2^8(1+2\alpha)}{\alpha\left({\alpha^2+2\alpha}\right)R_0^2 a}\left(\frac{2R}{R_0}\right)^{N-1}
 \|f_0\|_{{R_0}}\le 1
\end{eqnarray}
 any solution $(w(t), z(t))$ of the Hamilton equation of $H$ with initial datum  in $D(R)$ satisfying \eqref{reality at time zero} does not leave ${{D(R(1+\alpha))}}$ for all $t$ such that $|t|\le T_0$, where 
\begin{eqnarray}\label{stability time}
T_0&:=&\frac{1}{4}{(\alpha^2+2\alpha)}\min\left\{
\frac{a R_0^4\alpha^3}{2^{10}n(1+2\alpha) \|f\|^2_{{R_0}}}\left(\frac{R_0}{2R}\right)^{2N-2}\,,\ \right.\nonumber\\
&&
\left.\frac{\alpha R_0^2}{8\|f^{*}\|_{{R_0}}}\left(\frac{R_0}{2R}\right)^{M-1}\,,\ 
\frac{\alpha R^2}{2\|f_*\|_{{(1+2\alpha)R}}}
\right\}
\end{eqnarray}
where $f_{0}$, $f^*$ and $f_*$ are as in \eqref{fff}.
\end{theorem}

\begin{remark}\rm
The three terms in \eqref{stability time} come, respectively, from the contributions of the functions \eqref{fff}.
A possible choice for the parameter $M$, as justified by Theorem \ref{main simplified},  is $M=2N-1$. This corresponds to the case when the two former terms have comparable magnitude\footnote{Note however that the constants in front of them are very different, so, in applications  a finer evaluation may be  needed, see Theorem
\ref{main middle}.}  $1/R^{2N-2}$. Let us now analyze the third term 
in \eqref{stability time} whose magnitude is $1/R^{N-1}$, but 
as the parameter $a$ is small, the constant in front of him
is very large. Roughly, if $M$ is chosen as above, we may
re-write $T_0$ as
\begin{eqnarray}\label{minimum}
T_0=\frac{1}{4}{(\alpha^2+2\alpha)}\min\left\{
\frac{a R_0^4\alpha^3}{2^{10}n(1+2\alpha) \|f\|^2_{{R_0}}}\left(\frac{R_0}{2R}\right)^{2N-2}\,,\ \frac{\alpha R_0^2}{8\gamma(a)\|f\|_{{R_0}}}\left(\frac{R_0}{2R}\right)^{N-1}
\right\}
\end{eqnarray}
with $\gamma(a)\to 0$ as $a\to 0^+$. One expects that the minimum in Eq. \eqref{minimum} is reached by the first term, if $a$ is sufficiently small (note however that decreasing $a$ diminishes the first term, whence a too low value of $a$ might not be of help). In such a case, Theorem \ref{main simplified} provides a stability time $T_0$ which goes as $R^{-(2N-2)}$, hence 
improving the bound in \eqref{a priori bound}.
\end{remark}

Theorem \ref{main simplified} is a consequence of Theorem \ref{main}, stated
below. The latter also provides a potentially tighter bound for the action functions.

\begin{theorem}\label{main}
   Let $H$, as in \eqref{Hamiltonian}, be in Birkhoff normal form of order $N$ 
and analytic in $D(R_0)$. Let $M\in \mathbb N$, with    $M\ge N{+1}$, 
 $\alpha\in(0, {1/2}]$ and $a>0$. Then, for any  $0<R\le R_0/2$ verifying
   \eqref{smallness condition}, for any solution $(w(t), z(t))$ of the Hamilton equation such that $(w(t), z(t))\in {{D(R(1+\alpha))}}$ for all $t$ with $|t|\le T$, the following bound holds
\begin{eqnarray}\label{new bound}
\max_{j}\big|\big(I_j(T)-I_j(0)\big){- \big(\psi_{0j}(T)-\psi_{0j}(0)\big)\big|}&\le&\left(
\frac{64n(1+2\alpha)}{aR^2\alpha^3}\|f_{0}\|_{(1+2\alpha)R}\|f\|_{(1+2\alpha)R}+ \frac{2}{\alpha}\|f^*\|_{(1+2\alpha)R}\right.\nonumber\\
&+&\left.\frac{2}{\alpha}\|f_{*}\|_{(1+2\alpha)R}
\right)T
\end{eqnarray}
where $f_{0}$, $f^{*}$ and $f_*$ are as in \eqref{fff}, while $\psi_0(t)$ is given by
\begin{eqnarray}\label{psi(z,w)}
\psi_0(t):=\sum_{(h, k)\in{\cal I}_{ 0}}
\frac{c_{hk}(k-h)w(t)^h z(t)^k}{\Omega\cdot (h-k)}
\end{eqnarray}
and satisfies
\begin{eqnarray}\label{bound on psi0}
\sup_{|t|\le T} \max_{j}
|\psi_{0j}(t)|\le\frac{8(1+2\alpha)}{\alpha a}
\|f_0\|_{(1+2\alpha)R}  .
\end{eqnarray}

\end{theorem}

Additionally, we present a more technical version of Theorem \ref{main simplified}, which offers more flexibility for computer assisted approaches.

\begin{theorem}\label{main middle}
   Let $H$, as in \eqref{Hamiltonian}, in Birkhoff normal form of order $N$ 
and analytic in $D(R_0)$. Let $M\in \mathbb N$, with    $M\ge N{+1}$, 
 $\alpha\in(0, {1/2}]$ and $a>0$. Then, for any  $0<R\le R_0/2$ verifying
   \eqref{smallness condition}, and
\begin{eqnarray}\label{second smallness}\frac{2^4(1+2\alpha)}{\alpha(\alpha^2+2\alpha) a}
\frac{ \|f_0\|_{{(1+2\alpha)R}}}{R^2}< 1 ,
\end{eqnarray}
any solution $(w(t), z(t))$ of the Hamilton equation of $H$ with initial datum  in $D(R)$ satisfying condition \eqref{reality at time zero} does not leave ${{D(R(1+\alpha))}}$ for all $t$ such that $|t|\le T_1$, where 
\begin{eqnarray}\label{stability timeNEW}
T_1&:=& \left((\alpha^2+2\alpha)R^2-\frac{2^4(1+2\alpha)}{\alpha a}
\|f_0\|_{(1+2\alpha)R}\right) \left(
\frac{64n(1+2\alpha)}{aR^2\alpha^3}\|f_{0}\|_{(1+2\alpha)R}\|f\|_{(1+2\alpha)R}\right.\nonumber\\
&+& \frac{2}{\alpha}\|f^*\|_{(1+2\alpha)R}+\left.\frac{2}{\alpha}\|f_{*}\|_{(1+2\alpha)R}
\right)^{-1}
\end{eqnarray}
and $f_{0}$, $f^*$ and $f_{*}$ are as in \eqref{fff}.
\end{theorem}

\section{Proof of Theorems \ref{main simplified}, \ref{main} and \ref{main middle}}\label{The proofs}

The  proof of Theorems  \ref{main simplified}, \ref{main} and \ref{main middle}
exploits Eq. \eqref{integralwzII} for a selection
of multi-indices $(h,k)\in {\Bbb N}^n\times {\Bbb N}^n$. Indeed,
for all the solutions  $(w(\tau),z(\tau))$ which
belong to $D(R(1+\alpha))$ for all $\norm{\tau}\leq \norm{t}$, 
from the splitting \eqref{fff}, one can write
  \begin{eqnarray}\label{I(t)-I(0)}
 I_j(t)-I_j(0)&=&\int_0^t\{I_j, f\}d\tau\nonumber\\
    &=&
    \int_0^t\{I_j, f_0\}d\tau+ \int_0^t\{I_j, f_*\}d\tau+ \int_0^t\{I_j, f^*\}d\tau \,.   \end{eqnarray}
    The second and the third term will be treated in the traditional way, applying an upper bound of their modulus based on Cauchy inequalities. This will provide the two last terms in \eqref{new bound}.
    For what concerns the first term in \eqref{I(t)-I(0)},  instead of using Cauchy inequalities, we have the following
    \begin{lemma}\label{lemmafirstint}
 Let $H$, as in \eqref{Hamiltonian}, in Birkhoff normal form of order $N$ 
and analytic in $D(R_0)$. Let $M\in \mathbb N$, with    $M\ge N{+1}$, 
 $\alpha\in(0, {1/2}]$ and $a>0$. Then, for any  $0<R\le R_0/2$ verifying
  \eqref{smallness condition}, for any solution $(w(\tau), z(\tau))$ of the Hamilton equation satisfying condition \eqref{reality at time zero}
and $(w(\tau), z(\tau))\in {{D(R(1+\alpha))}}$ for all $|\tau|\le \norm{t}$,  we have
\begin{eqnarray}\label{firstintformula}
  \int_0^t\Big\{I_j, f_0\Big\}d\tau=\psi_j(w(t), z(t))-\psi_j(w(0), z(0))-\int_0^t\{\psi_j, f\}(w(\tau), z(\tau))d\tau
\end{eqnarray}
where \begin{eqnarray}\label{psiwz}
\psi(z,w):=\sum_{(h, k)\in{\cal I}_{ 0}}
\frac{c_{hk}(k-h)w^h z^k}{\Omega\cdot (h-k)}\,.
\end{eqnarray} 
\end{lemma}
    \begin{remark}\rm
      Notice that the vector function $\psi_0(t)$ in \eqref{psi(z,w)} satisfies  $\psi_0(t)=\psi(w(t), z(t))$.

\end{remark}

    \paragraph{\it Proof  of Lemma \ref{lemmafirstint}}  We  look
for an alternative  expression of the integrand function: 
    \begin{eqnarray} \label{Ijf0}
\Big\{I_j, f_0\Big\}=\ii\sum_{(h, k) \in{\cal I}_0}c_{hk}(k_j-h_j)w^h z^k
    \end{eqnarray}
    where we can use the small divisors hypothesis. Specifically, 
for all $(h, k) \in{\cal I}_0$ and all $\norm{\tau}\leq \norm{t}$ we have:
    \begin{equation}\label{smalldiv}
      \norm{\Omega(\ii w(\tau)z(\tau))\cdot(h-k)} \geq {a\over 2} .
    \end{equation}
   In fact, for all $(z,w)\in D(R(1+3\alpha/2))$ we have the estimate:
    \begin{eqnarray}|\Omega_i(\ii wz)-\Omega_{i}(0)|=
|\Omega_i(\ii wz)-\Omega_{0i}|=
      \left| \sum_{j=1}^n\partial_{I_j}\Omega_i(I^*)(\ii w_j z_j)\right|\le{L R^2(1+2\alpha)^2}\end{eqnarray}
with $L$ as in \eqref{L} and a suitable $I^*$ in the segment joining $I=0$ with
$I=iwz$. Since for all  $\norm{\tau}\leq \norm{t}$
we have $(w(\tau),z(\tau))\in D(R(1+\alpha))\subseteq  D(R(1+3\alpha/2))$, we
obtain:
 \begin{eqnarray} 
   \norm{\Omega(\ii w(\tau)z(\tau))\cdot(h-k)} &\geq& \norm{\Omega_0\cdot(h-k)} -
   \norm{(\Omega(\ii w(\tau)z(\tau))-\Omega_0)\cdot(h-k)}\cr
   &\geq& a -
   \sum_i \norm{   \Omega_i(\ii w(\tau)z(\tau))-\Omega_{0i}}  \norm{h_i-k_i}\cr
     &\geq& a - {L R^2(1+2\alpha)^2}\norm{h-k} \geq {a\over 2}
     \end{eqnarray}
 having used $\norm{h-k}\leq \norm{h}+\norm{k}\leq M$ and
 inequality \eqref{smallness condition}.

Then,  along any solution $(w(\tau), z(\tau))$, we eliminate $w^h z^k$ from \eqref{Ijf0} using the  identity
\begin{eqnarray}\label{whzk}
\frac{d}{dt} w^h z^k=\ii\Omega\cdot(h-k)w^h z^k+\Big\{w^h z^k, f\Big\}\end{eqnarray}
and we obtain (the division by $\ii\Omega(\ii w(\tau)z(\tau))\cdot(h-k)$
is allowed for all $ (h, k) \in{\cal I}_0$ because of \eqref{smalldiv})
 \begin{eqnarray}\label{developmentrhs}
\Big\{I_j, f_0\Big\}&=&\sum_{(h, k)\in{\cal I}_{ 0}}\frac{c_{hk}(k_j-h_j)}{\Omega\cdot (h-k)}\left(
\frac{d}{dt} w^h z^k-\Big\{w^h z^k, f\Big\}  
\right)  .\end{eqnarray}
 By integrating the functions at both sides of Eq. \eqref{developmentrhs} from $0$ to $t$,
 \begin{eqnarray}
\int_0^t \Big\{I_j, f_0\Big\}d\tau =
\sum_{(h, k)\in{\cal I}_{ 0}}\int_0^t \frac{c_{hk}(k_j-h_j)}
{\Omega\cdot (h-k)}\frac{d}{dt} w^h z^k d\tau
-\sum_{(h, k)\in{\cal I}_{ 0}}\int_0^t \frac{c_{hk}(k_j-h_j)}
{\Omega\cdot (h-k)} \{w^h z^k, f \}  d\tau  ,
\end{eqnarray}
by integrating by parts the first term at the right--hand side
 \begin{eqnarray}
   \int_0^t \Big\{I_j, f_0\Big\}d\tau &=&   
\sum_{(h, k)\in{\cal I}_{ 0}}\frac{c_{hk}(k_j-h_j)}
{\Omega\cdot (h-k)} w^h z^k  \Big \vert_0^t
-\sum_{(h, k)\in{\cal I}_{ 0}}\int_0^t w^h z^k 
\frac{d}{dt}
\frac{c_{hk}(k_j-h_j)}
{\Omega\cdot (h-k)} d\tau \cr
& & -\sum_{(h, k)\in{\cal I}_{ 0}}\int_0^t \frac{c_{hk}(k_j-h_j)}
{\Omega\cdot (h-k)} \{w^h z^k, f \}  d\tau
\end{eqnarray}
 and using:
 $$
 \frac{d}{dt}
\frac{c_{hk}(k_j-h_j)}
{\Omega\cdot (h-k)} = \left \{ \frac{c_{hk}(k_j-h_j)}
{\Omega\cdot (h-k)},Z+f \right \} = \left \{ \frac{c_{hk}(k_j-h_j)}
{\Omega\cdot (h-k)},f \right \}
$$
we obtain
 \begin{eqnarray}
   \int_0^t \Big\{I_j, f_0\Big\}d\tau &=&   
\sum_{(h, k)\in{\cal I}_{ 0}}\frac{c_{hk}(k_j-h_j)}
{\Omega\cdot (h-k)} w^h z^k  \Big \vert_0^t
-\sum_{(h, k)\in{\cal I}_{ 0}}\int_0^t w^h z^k 
\left \{ \frac{c_{hk}(k_j-h_j)}
{\Omega\cdot (h-k)},f \right \}
d\tau \cr
& & -\sum_{(h, k)\in{\cal I}_{ 0}}\int_0^t \frac{c_{hk}(k_j-h_j)}
{\Omega\cdot (h-k)} \{w^h z^k, f \}  d\tau  .
\end{eqnarray}
 Finally, using Leibniz rule to simplify the second and third terms
 in the right--hand side we obtain \eqref{firstintformula}. $\quad\square$

\paragraph{\it Some preliminary lemmas.} To prove Theorem \ref{main}, we first
provide some  preliminary lemmas, which, even though being standard, are nevertheless  necessary to provide the exact value of the constants appearing in \eqref{new bound}.

\vskip.1in
\noindent
In addition to the norm in \eqref{Pnorm}, we also define
\begin{eqnarray}\label{inftynorm}|f|_R:={{\max}_{{{D(R)}}}|f|}\,.
\end{eqnarray}
First of all, we observe that the norms in \eqref{Pnorm}, \eqref{inftynorm} satisfy the trivial bound
\begin{lemma}
\label{supbound} $\displaystyle| f|_R\le\|f\|_R
$.
 \end{lemma}
In lemmas \ref{derivatives}--\ref{bound on actions}  below we assume that
$f$, $g$ are analytic functions on ${{D(R_0)}}$,  $0<R\le R_0$, $\alpha$, $\beta\ge 0$, $0<\alpha+\beta\le\frac{R_0-R}{R}$.
 \begin{lemma}\label{derivatives} 
For both norms in \eqref{Pnorm}, \eqref{inftynorm}, \begin{eqnarray}\|\partial_{\zeta_i} f\|_{(1+\alpha)R}\le \frac{\|f\|_{(1+\alpha+\beta)R}}{\beta R}\end{eqnarray}
  \end{lemma}
\paragraph{\it Proof of Lemma \ref{derivatives} }
  The  bound  is classical for $|\cdot|_R$, so it will not be discussed. As for 
$\|\cdot\|$
\begin{eqnarray}
\|\partial_{w_i} f\|_{(1+\alpha)R}&=&\sum_{h, k
 }h_i|c_{hk}|R^{|h|+|k|-1}(1+\alpha)^{|h|+|k|-1}\nonumber\\
 &\le& \frac{1}{\beta R}\sum_{h, k
 }|c_{hk}|R^{|h|+|k|}(1+\alpha+\beta)^{|h|+|k|}\nonumber\\
 &=&\frac{\|f\|_{(1+\alpha+\beta)R}}{\beta R}\end{eqnarray}
having used
$1+na\le (1+a)^n$. The case $\zeta_i=z_i$ is the same. $\quad \square$
\begin{lemma}\label{Poisson brackets} For both norms in \eqref{Pnorm}, \eqref{inftynorm},
\begin{eqnarray}\|\{f, g\}\|_{(1+\alpha)R}\le \frac{2n}{\beta^2 R^2}\|f\|_{(1+\alpha+\beta)R}\|g\|_{(1+\alpha+\beta)R}\end{eqnarray}
\end{lemma}
\paragraph{\it Proof of Lemma \ref{Poisson brackets}} It follows from Lemma \ref{derivatives}. $\quad\square$
\begin{lemma}\label{bound on actions}
For both norms in \eqref{Pnorm}, \eqref{inftynorm},
\begin{eqnarray}
\|\{I_j, f\}\|_{(1+\alpha)R}\le\frac{2(1+\alpha)}{\beta}\|f\|_{(1+\alpha+\beta)R}
\end{eqnarray}
\end{lemma}
\paragraph{\it Proof of Lemma \ref{bound on actions}} We have
\begin{eqnarray}
\{I_j, f\}=\ii \left(z_j \partial_{z_j} f-w_j \partial_{w_j} f\right)\end{eqnarray}
Thus, by Lemma \ref{derivatives},
\begin{eqnarray}
\|\{I_j, f\}\|_{(1+\alpha)R}&\le& \|z_j\|_{(1+\alpha)R}\| \partial_{z_j} f\|_{(1+\alpha)R}+\|w_j \|_{(1+\alpha)R}\|\partial_{w_j} f\|_{(1+\alpha)R}\nonumber\\
&\le&\frac{2(1+\alpha)}{\beta}\|f\|_{(1+\alpha+\beta)R}\,.\qquad \square \end{eqnarray}

\vskip.1in
\noindent
\paragraph{\it Proof  of Theorem \ref{main}} Let us consider any
solution $(w(t), z(t))$ of the Hamilton equation  such that $(w(t), z(t))\in {{D(R(1+\alpha))}}$ for all $|t|\le T$. We 
represent the first integral on the right hand side of Eq. (\ref{I(t)-I(0)})
for $t=T$:
  \begin{eqnarray}\label{represintegrals}
 I_j(T)-I_j(0) = 
 \int_0^T\{I_j, f_0\}d\tau+ \int_0^T\{I_j, f_*\}d\tau+ \int_0^T\{I_j, f^*\}d\tau \,.   \end{eqnarray}
  using Lemma  \ref{lemmafirstint} (whose hypotheses are satisfied with $t=T$),
  and obtain (see Eq. \eqref{firstintformula}):
\begin{eqnarray}\label{firstintformulaT}
  \int_0^T\Big\{I_j, f_0\Big\}d\tau=\psi_{0j}(T)-\psi_{0j}(0)-\int_0^T\{\psi_j, f\}(w(\tau), z(\tau))d\tau .
\end{eqnarray}

\noindent
Therefore, using Cauchy inequalities to estimate the second and third
integrals in the right-hand side of Eq. \eqref{represintegrals} and the last integral  in the right-hand side of Eq. \eqref{firstintformulaT}, we have the following bound
\begin{eqnarray}\label{three terms}\big|\big(I_j(T)-I_j(0)\big)-\big(\psi_{0j}(T)-\psi_{0j}(0)\big)\big|\le(|\{\psi_j, f\}|_{(1+\alpha)R}+|\{I_j, f_*\}|_{(1+\alpha)R}+|\{I_j, f^*\}|_{(1+\alpha)R})T\,.
\end{eqnarray}
We now proceed to provide upper bounds on the norms of the terms at right hand side. For the second and the third term, we use Lemma \ref{bound on actions} with
 $\beta=\alpha$.
We obtain
\begin{eqnarray}\label{second and third}
|\{I_j, f_*\}|_{(1+\alpha)R}\le\frac{2(1+\alpha)}{\alpha}|f_*|_{(1+2\alpha)R}\,,\qquad |\{I_j, f^*\}|_{(1+\alpha)R}\le\frac{2(1+\alpha)}{\alpha}|f^*|_{(1+2\alpha)R}\,.
\end{eqnarray}
To bound the first term in \eqref{three terms}, we first use  Lemma \ref{Poisson brackets} with  $\beta=\frac{\alpha}{2}$. This gives
\begin{eqnarray}\label{psif}|\{\psi_j, f\}|_{(1+\alpha)R}\le\frac{8n}{\alpha^2R^2}|\psi_j|_{(1+\frac{3}{2}\alpha)R}|f|_{(1+\frac{3}{2}\alpha)R} .
\end{eqnarray}
Now, since for all $(z,w)\in {{D(R(1+3\alpha/2))}}$
and all $(h, k)\in{\cal I}_0$
we have (see the Proof of Lemma \ref{lemmafirstint})
\begin{eqnarray}\label{lower bound}
  |\Omega\cdot (h-k)|
\ge \frac{{a}}{2} ,
\end{eqnarray}
we get
\begin{eqnarray}\label{psi}|\psi_j|_{(1+\frac{3}{2}\alpha)R}\le \frac{2}{a}\|g_{0j}\|_{(1+\frac{3}{2}\alpha)R}\end{eqnarray}
where
\begin{eqnarray}
\label{gstar} g_{0j}:=\sum_{(h, k)\in {\cal I}_0}c_{hk}(h_j-k_j)w^h z^k=\ii \{I_j, f_0\}\,.
\end{eqnarray}
Using Lemma \ref{bound on actions} again,  with $\alpha$ replaced by $\frac{3}{2}\alpha$ and $\beta=\frac{\alpha}{2}$, we have
\begin{eqnarray}\label{g0}
\|g_{0j}\|_{(1+\frac{3}{2}\alpha)R}=\|\{I_j, f_0\}\|_{(1+\frac{3}{2}\alpha)R}\le 4\frac{1+2\alpha}{\alpha}\|f_0\|_{(1+2\alpha)R}
\end{eqnarray}
having replaced, at right hand side, $\frac{3}{2}\alpha$ with $2\alpha$ by simplicity.
The bounds in \eqref{psi} and \eqref{g0} imply
\begin{eqnarray}
|\psi_j|_{(1+\alpha)R}\le\frac{8(1+2\alpha)}{\alpha a}
\|f_0\|_{(1+2\alpha)R} .
\end{eqnarray}
As $(w(t), z(t))\in D(R(1+\alpha))$ for all $|t|\le T$, inequality \eqref{bound on psi0} follows for $\psi_0(t)=\psi(w(t), z(t))$.
Moreover, 
\eqref{psif}, \eqref{psi} and \eqref{g0} lead to
\begin{eqnarray}\label{quadratic}|\{\psi_j, f\}|_{(1+\alpha)R}\le\frac{64n(1+2\alpha)}{\alpha^3R^2a}
|f|_{(1+2\alpha)R}\|f_0\|_{(1+2\alpha)R}
\end{eqnarray}
Finally, \eqref{three terms}, \eqref{second and third} and \eqref{quadratic} imply
\begin{eqnarray}
\big|\big(I_j(T)-I_j(0)\big){- \big(\psi_{0j}(T)-\psi_{0j}(0)\big)\big|}&\le&\left(
\frac{64n(1+2\alpha)}{aR^2\alpha^3}\|f_{0}\|_{(1+2\alpha)R}|f|_{(1+2\alpha)R}+ \frac{2}{\alpha}|f^*|_{(1+2\alpha)R}\right.\nonumber\\
&+&\left.\frac{2}{\alpha}|f_{*}|_{(1+2\alpha)R}
\right)T
\end{eqnarray}
which in turn implies \eqref{new bound}, using Lemma \ref{supbound}. $\quad \square$

\paragraph{\it Proof  of Theorem \ref{main middle}} We first notice that, 
from Eq. \eqref{second smallness}, the time $T_1$ defined in
Eq. \eqref{stability timeNEW} is strictly positive. 
Let us   consider a solution  $(w(t), z(t))$ of the Hamilton equations of $H$  with initial datum in ${{ D(R)}}$, satisfying condition \eqref{reality at time zero}. Let $T_{\rm ex}$ the minimum value of $t>0$ such that
 \begin{eqnarray}
 |w_j(T_{\rm ex})|{=|z_j(T_{\rm ex})|}=R(1+\alpha)
 \end{eqnarray} for some $j=1$, $\ldots$, $n$,  as well as
\begin{eqnarray}\label{boundaryNEW}
I_j(T_{\rm ex})=\ii w_j(T_{\rm ex})z_j(T_{\rm ex})=|w_j(T_{\rm ex})|^2=R^2(1+\alpha)^2\,.\end{eqnarray}
{By the definition of $T_{\rm ex}$, }, for all $t$ with $ |t|\le T_{\rm ex}$, the solution $(w(t), z(t))$ stays in ${{D((1+\alpha)R)}}$, hence,  Theorem \ref{main} applies with $T=T_{\rm ex}$.

\noindent
Using  \eqref{new bound}, \eqref{bound on psi0} and the triangular inequality we have
\begin{eqnarray}
\max_j \big|\big(I_j(T_{{\rm ex}})-I_j(0)\big)\big|
&\le&\frac{16(1+2\alpha)}{\alpha a}
 \|f_0\|_{{(1+2\alpha)R}}+
 \left(
\frac{64n(1+2\alpha)}{aR^2\alpha^3}\|f_{0}\|_{(1+2\alpha)R}\|f\|_{(1+2\alpha)R}\right.\nonumber\\
&+& \frac{2}{\alpha}\|f^*\|_{(1+2\alpha)R}+\left.\frac{2}{\alpha}\|f_{*}\|_{(1+2\alpha)R}
\right)T_{\rm ex}\end{eqnarray}
Assume now, by contradiction, that 
$T_{\rm ex}$ is {strictly smaller} than the number $T_1$  in \eqref{stability timeNEW}. 
Then,
\begin{eqnarray}\label{strict inequality}
|I_j(T_{{\rm ex}})|&\le& |I_j(0)|+|I_j(T_{{\rm ex}})-I_j(0)|< R^2+
{(\alpha^2+2\alpha)}
R^2=R^2\left(1+{\alpha}\right)^2
\end{eqnarray}
for all  $j=1$, $\ldots$, $n$.
The inequality in \eqref{strict inequality}
contradicts \eqref{boundaryNEW}. $\quad \square$

\paragraph{\it Proof of Theorem \ref{main simplified}}

From the assumption $\alpha\in (0, \frac{1}{2}]$, one has $R(1+2\alpha)\le 2R$. 
 Moreover,
\begin{eqnarray}\label{norm of f}
\|f\|_{R(1+2\alpha)}&\le&\|f\|_{ 2R}=
\sum_{|h|+|k|\ge N+1} |c_{hk}| \left( 2R\right)^{|h|+|k|}\nonumber\\
&=&
\sum_{|h|+|k|\ge N+1} |c_{hk}| \left(\frac{ 2R}{R_0}\right)^{|h|+|k|}
\left(R_0\right)^{|h|+|k|}
\nonumber\\
&\le&\left(\frac{2R}{R_0}\right)^{N+1}
 \|f\|_{{R_0}}
\end{eqnarray}
Similarly,
\begin{eqnarray}\label{norm of f0}
&&\|f_0\|_{R(1+2\alpha)}\le\left(\frac{2R}{R_0}\right)^{N+1} \|f_0\|_{{R_0}}\cr
&&\|f^*\|_{R(1+2\alpha)}\le\left(\frac{2R}{R_0}\right)^{M+1} \|f^*\|_{{R_0}}\,.
\end{eqnarray}
The inequalities \eqref{norm of f0} and \eqref{smallness condition NEW} imply \eqref{second smallness}. Then Theorem \ref{main middle} applies. Moreover, by \eqref{norm of f}, \eqref{norm of f0}  the number $T_1$ in \eqref{stability timeNEW} is larger or equal than the number $T_0$ in \eqref{stability time}. $\qquad\square$

 \section{Rigorous numerical application to a 3 d.o.f. case}\label{application}

 We apply the results of Section \ref{thestatements} to study the stability of
 the equilibrium point of a three-degrees of freedom system, obtained from
 the well known Fermi-Pasta-Ulam $\alpha$-model Hamiltonian 
 for 5 interacting particles (the two at the boundary are fixed)\footnote{
Although the equilibrium point in this model is Lyapunov stable due to the positivity of all harmonic oscillation frequencies, we employ it to compare stability times derived from classical estimates with those obtained using this new technique. Furthermore, these stability estimates serve to establish refined upper bounds on actions variations, despite the inherent Lyapunov stability of the system.}. All
 computations are performed using the algebraic manipulator {\grm Qr'onoc}~\cite{giosan} and made rigorous by the use of interval arithmetics.

 \paragraph{\it The  3 d.o.f. FPU $\alpha$-model and its Birkhoff normal forms.}  By denoting with $(y_j,x_j)$
 the momentum and  the configuration of the $j$-th particle, $j=0,\ldots ,4$,
and by constraining $x_0 = x_4 = y_0 = y_4 = 0$,
we have the 3-degrees of freedom FPU $\alpha$-model Hamiltonian:
\begin{equation}\label{hamFPU0}
H_0(y,x) = {1\over 2}\sum_{j=0}^{3}\left[y_j^2+(x_{j+1}-x_j)^2\right] + \frac{\tilde \alpha}{3} \sum_{j=0}^3(x_{j+1}-x_j)^3 
\end{equation}
defined for $(y,x):=(y_1,y_2,y_3,x_1,x_2,x_3)\in {\Bbb R}^{6}$. 
The origin $(y,x)=(0,\ldots, 0)$ is an elliptic equilibrium point,  and the canonical change of coordinates 
$$
x_\ell = \frac{1}{\sqrt{2}} \sum_{j=1}^3 \frac{q_j}{\sqrt{\Omega_{0j}}} \sin\left(\frac{j\ell \pi}{4}\right), \quad {\rm and } \quad
y_\ell = \frac{1}{\sqrt{2}} \sum_{j=1}^3 p_j\sqrt{\Omega_{0j}} \sin\left(\frac{j\ell \pi}{4}\right)\ \ ,\ \ \ell=1,2,3
$$
with 
\begin{equation}\label{Omega0FPU}
\Omega_{0j} = 2 \sin\left(\frac{j \pi}{8}\right)
\end{equation}
conjugates Hamiltonian (\ref{hamFPU0}) to
\begin{equation}\label{FPU1}
H(p,q)= \sum_{j=1}^3  \Omega_{0j} \frac{p_j^2+q_j^2}{2} + \frac{\tilde \alpha}{3} H_3(p,q)
\end{equation}
where $H_3$ is an homogeneous polynomial of degree 3; as in classical studies, we take $\tilde \alpha=1/4$. 

We introduce Poincar\'e complex variables as in equation~\eqref{wz}, and
 conjugate
Hamiltonian (\ref{FPU1}) to its Birkhoff normal forms  \eqref{Hamiltonian} of order\footnote{
The construction is rather standard in perturbation theory and it is done by the means of near-the-identity canonical change of coordinates through Lie series . We defer to~\cite{carloc2020} for a detailed description of the algorithm and how to perform rigorous numerics.} $N$ :
\begin{equation}
\label{normalizedhh}
{\cal H}(z,w) = Z(\ii zw) + f(z,w),
\end{equation}
where the remainder is split in $f= f_0+f_*+f^*$, according to the definitions
given in~\eqref{fff} and~\eqref{indices}. Hence, the terms $f_0$ and $f_*$ are polynomial terms of maximal degree $M$ which are explicitly computed with an algebraic manipulator, while $f^*$ is not computed but its norm is bounded with \emph{a priori} analytical estimates. The distinction between the terms in $f_0$ and $f_*$ depends on the value of a parameter $a$, which will be conveniently chosen later.

 \paragraph{\it The  classical estimate of the stability time.}  We aim to compare the stability times computed from formula~\eqref{stability timeNEW} 
of Theorem~\ref{main middle} with the stability times which are computed
using standard {\it a priori} estimates (hereafter denoted as {\it classical}
estimates). Let us review how the classical estimates are calculated. For now, let us collect in $f_0$ all the explicit terms appearing in the remainder (and take $f_*=0$). 
For what concerns the estimate of $f^*$ (whose terms are not
explicitly computed), there exists a positive real constant  $b_N$ such that
\begin{equation}
\label{tail_est}
|f^*|_\rho\le \|f^*\|_\rho\le \sum_{j=M+1}^\infty (b_N \rho)^{j} = \frac{(b_N \rho)^{M+1}}{1-b_N \rho},
\end{equation}
which we compute following~\cite{carloc2020}, Proposition 2. Since
the series appearing in the right--hand side of inequality (\ref{tail_est}) is convergent for $|b_N \rho|<1$,
the value of $b_N$ poses a constraint on the radius of the domain in which
estimate (\ref{tail_est}) is valid. From (\ref{tail_est}), we then
recover the inequality
\begin{equation}
| f|_\rho \le  \|f_0+f_*\|_\rho+ \|f^*\|_\rho
\end{equation}
since $\|f_0+f_*\|_\rho$ is explicitly computed (the coefficients of all
the monomials in the Birkhoff normal form with degree less than or equal to $M$ are known). Finally, since the time variation of the actions satisfies
\begin{equation}
\label{classicavaraz}
|I_j(T)-I_j(0)|_\rho \le \int_0^T |\{I_j,f\}|_\rho dt \le T |\{I_j,f\}|_\rho\,,
\end{equation}
using Lemma~\ref{bound on actions} and taking $\rho =(1+\alpha)R$, we
first obtain
\begin{equation}
\label{deltaIclassic}
|\{I_j,f\}|_{(1+\alpha)R} \le \|\{I_j,f\}\|_{(1+\alpha)R}\le  \frac{2(1+\alpha)}{\alpha}\|f\|_{(1+2\alpha)R},
\end{equation}
and then we define a lower bound for the time needed to escape from $D(R(1+\alpha))$ as\footnote{We impose $|I_j(T)|\le |I_j(0)|+|I_j(T)-I_j(0)| \le R^2(1+\alpha)^2$. Therefore, since $|I_j(0)|<R^2$, we choose $T$ so as to satisfy the inequality $|I_j(T)-I_j(0)|  \le R^2\alpha(2+\alpha)$.}
\begin{equation}
\label{tempost}
T_{c}:=\frac{R^2\alpha^2(2+\alpha)}{2(1+\alpha)\|f\|_{(1+2\alpha)R}}.
\end{equation}
The  stability time $T_c$ will be hereafter identified  as the classical estimate. 

\paragraph{\it Comparison of the stability times.} Given $R,\alpha,N, M$ and  $R_0{\leq 1/b_N}$ the stability time $T_c$ 
is compared with the stability time $T_1$ computed using formula~\eqref{stability timeNEW} of Theorem~\ref{main middle}, with conveniently chosen parameter
$a$ satisfying conditions \eqref{smallness condition} and \eqref{second smallness}. For any given $R$, we compute $T_1$ with parameter $a$ chosen as
\begin{equation}
\label{a}
a=\max \left\{2LM(1+2\alpha)^2 R^2, \min_{N<|k|+|h|\le M,{k \ne h}} |(k-h)\cdot \Omega_0| \right\},
\end{equation}
so that condition~\eqref{smallness condition} is
satisfied, and we directly check if condition~\eqref{second smallness} is satisfied as well (typically, the condition is satisfied in view of the smallness of $\|f_0\|_{(1+2\alpha)R}$). In Table~\ref{tab1} we report for different values of
$R$ the stability times $T_c$ and $T_1$ computed using formulae~\eqref{tempost} and~\eqref{stability timeNEW} for $\alpha=1/5$, $N=9$, two diffent values of $M$: $M=2N= 18$ and $M=3N =27$, and:
\begin{equation}
\label{R0}
R_0:= \frac{1}{b_N}.
\end{equation}

\begin{remark}\rm Let us justify the choice of the parameters $a,M,R_0$.  
  \begin{itemize}
\item[-]  When $R$ is particularly small,
  the value $a_0:=2LM(1+2\alpha)^2 R^2$ is so small that typically there are no divisors smaller than it, and consequently the choice $a=a_0$
  would imply $f_* =0$. However, since $T_1$ is singular for $a=0$, 
  such a small value of $a$ has an impact in the estimate of the stability time $T_1$. Therefore, we find convenient to increase $a$ to the value of the minimum divisor that we have among the terms of the remainder which are explicitly computed.
\item[-]  The choice of $M$ was, in first analysis, motivated by the asymptotic analysis in the remark \ref{remarkM=2N}, which emphasizes how the stability time goes as the minimum between a term of order $R^{2-2N}$ and one of order $R^{1-M}$. Therefore, one would not expect to get any improvement by taking $M> 2N-1$. However, from a computational point of view, we must consider that, while the norm $\|f_0\|_{(1+2\alpha) R}$ is computed exactly (in the limit of interval arithmetics), $\|f^*\|_{(1+2\alpha R)}$ is only \emph{estimated} using Equation~\eqref{tail_est}. It is clear from Equation~\eqref{tail_est} that the greater $M$, the smaller will be the estimate for $\|f^*\|_{(1+2\alpha) R}$, because the number of terms that that are explicitly computed increases. Due to this computational issue, by taking $M=2N$ we typically found that, in Formula~\eqref{stability timeNEW}, the terms depending on $\|f^*\|_{(1+2\alpha) R}$ is much greater than the one in $\|f_0\|_{(1+2\alpha) R}\| f\|_{(1+2\alpha) R}$. In order to better balance the two contributions, we tested a value of $M=3N$. The comparison between times in Table~\ref{tab1} shows how the stability time $T_1$ is much longer than before and it is greater than $T_{c}$ for higher values of $R$.
\item[-] The last values of $R$ reported in Table~\ref{tab1} are limit values for the method, with that choices of $N$ and $M$.  Indeed, we recall that, in order to get a convergent series for the estimate of $\|f^*\|_{R(1+2\alpha)}$ in~\eqref{tail_est}, the condition $b_N R(1+2\alpha)<1$ must be satisfied. Since $b_N$ is an increasing function of $N$, the bigger the $N$, the smaller will be the maximal radius $R$ for which we can apply the estimates.
  This leads to the definition of the limit value in (\ref{R0}), which
  is  reported below Table~\ref{tab1}.
 \end{itemize} 
\end{remark}

\begin{table}[t]
    \centering
        \begin{tabular}{c|ccc||ccc}
 \hline
 & & M=18  & & &  M=27& \\ 
       $R$ & $\log_{10}(T_{c})$& $a$ & $\log_{10}(T_{1)}$ &  $\log_{10}(T_{c})$ &$a$ & $\log_{10}(T_{1)}$\\
 \hline
5.000000e-04 & 23.825453 & 1.185291e-02 & 31.832795 & 23.825453 & 5.886225e-03 & 43.409313 \\
6.000000e-04 & 23.191935 & 1.185291e-02 & 30.486008 & 23.191935 & 5.886225e-03 & 42.142276 \\
7.200000e-04 & 22.558403 & 1.185291e-02 & 29.139080 & 22.558403 & 5.886225e-03 & 40.875212 \\
8.640000e-04 & 21.924854 & 1.185291e-02 & 27.791979 & 21.924854 & 5.886225e-03 & 39.608111 \\
1.036800e-03 & 21.291283 & 1.185291e-02 & 26.444672 & 21.291286 & 5.886225e-03 & 38.340946 \\
1.244160e-03 & 20.657676 & 1.185291e-02 & 25.097115 & 20.657695 & 5.886225e-03 & 37.073584 \\
1.492992e-03 & 20.023976 & 1.185291e-02 & 23.749255 & 20.024075 & 5.886225e-03 & 35.805238 \\
1.791590e-03 & 19.389912 & 1.185291e-02 & 22.401031 & 19.390420 & 5.886225e-03 & 34.531117 \\
2.149908e-03 & 18.754095 & 1.185291e-02 & 21.052364 & 18.756725 & 5.886225e-03 & 33.223580 \\
2.579890e-03 & 18.109490 & 1.185291e-02 & 19.703159 & 18.122981 & 5.886225e-03 & 31.765754 \\
3.095868e-03 & 17.423194 & 1.185291e-02 & 18.353300 & 17.489177 & 5.886225e-03 & 29.991204 \\
3.715042e-03 & 16.587015 & 1.185291e-02 & 17.002641 & 16.855302 & 5.886225e-03 & 27.998950 \\
4.458050e-03 & 15.483997 & 1.185291e-02 & 15.651002 & 16.221342 & 5.886225e-03 & 25.947221 \\
5.349660e-03 & 14.200764 & 1.185291e-02 & 14.298157 & 15.587279 & 5.886225e-03 & 23.883784 \\
6.419592e-03 & 12.861107 & 1.185291e-02 & 12.943817 & 14.953092 & 5.886225e-03 & 21.817139 \\
7.703511e-03 & 11.507760 & 1.185291e-02 & 11.587613 & 14.318755 & 5.886225e-03 & 19.748356 \\
9.244213e-03 & 10.149757 & 1.185291e-02 & 10.229065 & 13.684188 & 5.886225e-03 & 17.677185 \\
1.109306e-02 & 8.788331 & 1.185291e-02 & 8.867536 & 13.048054 & 5.886225e-03 & 15.603025 \\
1.331167e-02 & 7.422975 & 1.185291e-02 & 7.502160 & 12.375878 & 5.886225e-03 & 13.525019 \\
1.597400e-02 & 6.052553 & 1.185291e-02 & 6.131735 & 11.221854 & 5.886225e-03 & 11.441962 \\
1.916880e-02 & 4.675334 & 1.185291e-02 & 4.754515 & 9.267119 & 5.886225e-03 & 9.352112 \\
2.300256e-02 & 3.288678 & 1.185291e-02 & 3.367859 & 7.173441 & 5.886225e-03 & 7.252823 \\
2.760307e-02 & 1.888302 & 1.185291e-02 & 1.967483 & 5.060629 & 6.171355e-03 & 5.139816 \\
3.312369e-02 & 0.466602 & 1.185291e-02 & 0.545781 & 2.926305 & 8.886751e-03 & 3.005483 \\
3.974842e-02 & -0.991993 & 1.185291e-02 & -0.912824 & 0.755078 & 1.279692e-02 & 0.834248 \\
        \end{tabular}
        \caption{Stability times for different values of $R$ and $N=9$, $M=18,27$; $R_0 =1/b_9= 8.700956$e-02.}
\label{tab1}
\end{table}

\begin{figure}[!]
\centering
\includegraphics[scale=0.57]{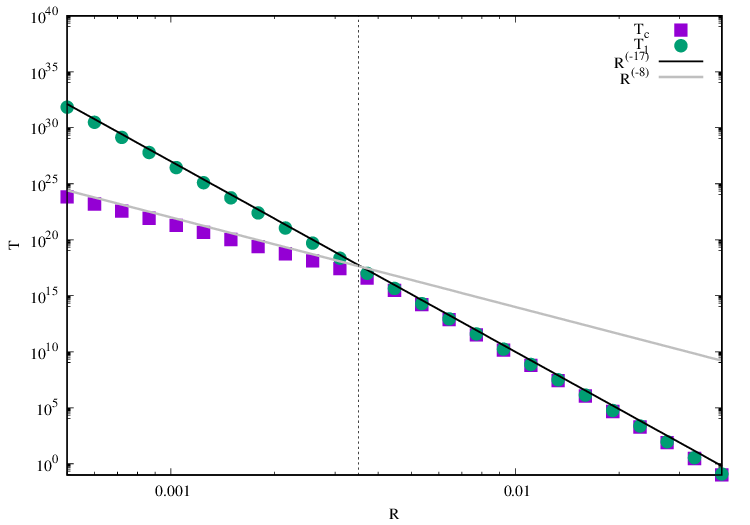}
\includegraphics[scale=0.57]{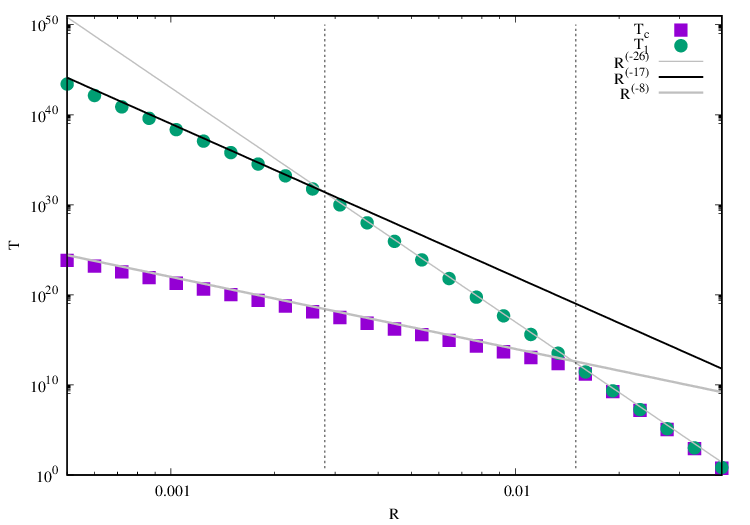}
\caption{Logscale plot of the function $T_1(R)$ and $T_c(R)$ with $N=9$ and $M=18$ (left), $M=27$ (right). The lines highlight the decrease rates: $T_{c}\simeq R^{1-N}$; $T_1\simeq\min\{c_1 R^{1-2N},c_2 R^{1-M}\}$. In the plot on the right one can identify two regimes for $T_1$: for small $R$, it goes as $R^{-17}$, while after some time term $R^{-26}$ becomes predominant. The two estimates come to coincide when $\|f^*\|\ge \|f_0+f_*\|$. }
\label{fig}
\end{figure}

  Let us know comment Table~\ref{tab1}, where, for each value of $R$,
  we report the values of $T_c,T_1,a$ computed for $M=18$ and $M=27$
  respectively.

\begin{itemize}
\item For small values of $a$ one has $\mathcal I_*= \emptyset$ and $f_* =0$. However, increasing $R$ also $a$ increases and there could be divisors smaller than $a$, that is, when the $\max$ in formula~\eqref{a} is obtained from the first term. For example, for $M=27$, $I_*$ is non-empty and $f_* \neq 0$ when $R \ge 2.760307$e-02 because we have the divisor
$\norm{(h-k)\cdot \Omega_0 =( -4,10,-6)\cdot \Omega_{0}} = 0.005886...$
which becomes smaller than $a=0.006175...$ 
Instead, for $M=18$ not only $a$ grows slower, but there could be less divisors, so that,  $f_*$ is always zero. 

\item The values of $R$ displayed  in the Table have been increased at the constant rate of $1.2$, in order to better appreciate the trend. The decrease of the stability time with respect to the radius $R$ is emphasized by the logscale plot in Figure~\ref{fig}, where one can appreciate different regimes for $T_{c}$ and $T_1$. In particular, we observe how $T_{c}\simeq R^{1-N}$, as we expect theoretically, up to a limit value of $R$ (emphasized by the vertical line on the right), from which it started to decrease as $R^{1-M}$; $T_1$ initially goes as $R^{2-2N}$ and after a specific value of $R$ (the vertical line on the left) it goes as $R^{1-M}$. We can observe also that the two stability times $T_1$ and $T_{c}$ come to coincide when the estimate of the norm of $f^*$ becomes the main contributor to the estimate of the norm of the total remainder $f$. In fact, we recall that the improvement in our method is made on the explicit part $f_0$, therefore the estimates come to coincide when $f^*$ becomes predominant. This limit value of $R$ for which the two estimates coincide could be increased by increasing $M$ or improving the estimate of $\|f^*\|_{(1+2\alpha)R}$.
  \end{itemize}

\begin{figure}[!]
\centering
\includegraphics[scale=0.6]{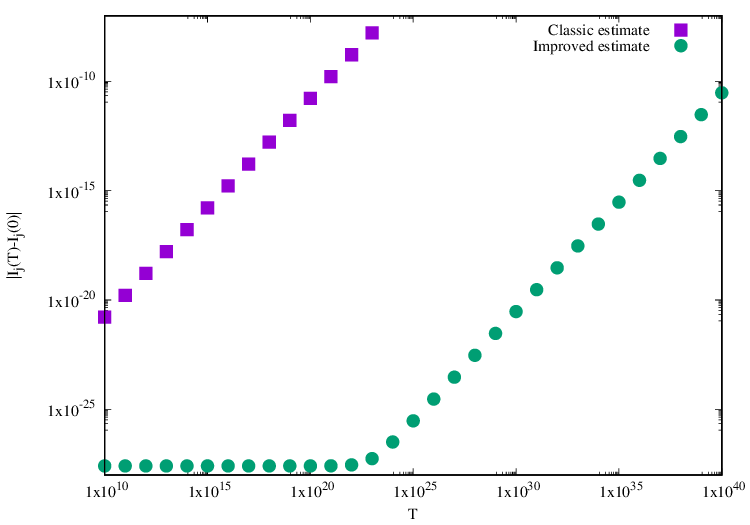}
\caption{Logscale plot of the estimated actions variation with the classical (violet squares) and new (green circles) estimate, for different times $T$. The initial condition satisfies $I_j(0)<R^2$ with $R=5.$e-4.}
\label{fig2}
\end{figure}

Finally, we reported in Figure~\ref{fig2} the estimate in the variation of the actions $ |I_j(T)-I_j(0)|$, taking $I_j(0)<R^2$ with $R=5.$e-4 for different times $T$, using the classical estimate~\eqref{classicavaraz} and the estimate given in Theorem~\ref{main}, formula~\eqref{new bound}. We can observe how the same variations of the actions occurs in much longer times than with classical estimates. We can also recognize here the different regimes for the new estimates. Indeed, in formula~\eqref{new bound}, the time $T$ multiplies terms that are much smaller than in the standard case (quadratic in the remainder or of order $\|f^*\|$, instead of all $\|f\|$). This means that for small values of $T$, it is predominant the stationary phase term $|\psi_j(T)-\psi_j(0)|$.

\begin{remark}\rm
In \cite{carloc2020} alternative strategies for the computation of the stability time are used, in particular using two cut-offs
$M_1< M_2$ larger than $N$, rather than a single cut-off $M$. The comparison
discussed in this paper refers to the simpler strategy. However, we emphasize that the refined strategy used in \cite{carloc2020} could only
improve the estimate of $\|f^*\|_{(1+2\alpha)R}$, which contributes equally
to both the classical and the improved estimates of the stability time.
Therefore, including the method described in \cite{carloc2020} could improve the
stability time, but not the comparison between the two methodologies.
\end{remark}

\section{On the size of $\|{ f}_*\|_{2R}$}\label{An upper bound on the resonant term}
Arguing as done in the proof of Theorem \ref{main simplified} to obtain the formulae \eqref{norm of f}--\eqref{norm of f0}, one can infer that the  term $\|{ f}_*\|_{2R}$ appearing in the formula \eqref{stability time} is of size $C_N(a)\left(\frac{2R}{R_0}\right)^{N+1}$, where $C_N(a)$ is a coefficient related not
only to $N$, but also to the number of couples $(h, k)$ in ${\cal I}_*$. 
So, even though the worst scenario assigns to such term a size $\sim \left(\frac{2R}{R_0}\right)^{N+1}$, on the other side, it is reasonable to expect that the coefficient $C_N(a)$ is so small to make such term of negligible size. The purpose of this Section is to provide some rigorous assertion on this question {for the
  case $n=3$}, which, for what said, is of central importance in the present paper.
\vskip 0.4 cm

\begin{figure}[!]
\centering
\includegraphics[scale=0.57]{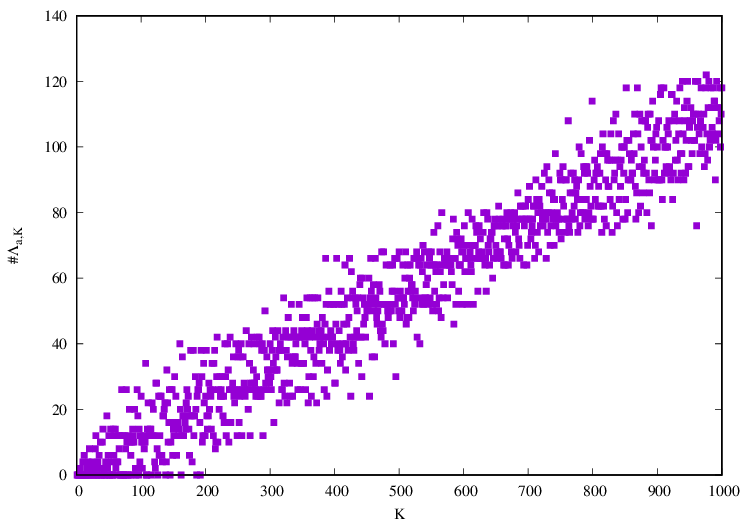}
\includegraphics[scale=0.57]{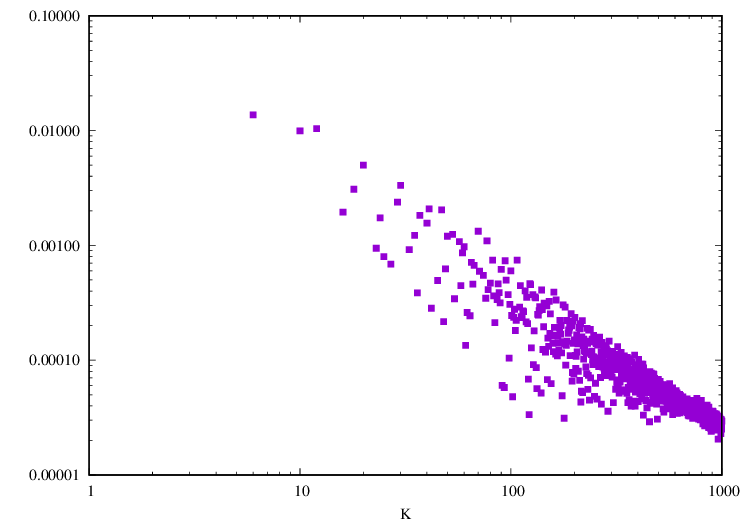}
\caption{Representation of the  cardinality $\#\Lambda_{a, K}{(\Omega_0)}$
  for the vector $\Omega_0\in{\Bbb R}^3$ defined in Section 4, Eq. (\ref{Omega0FPU}) (left panel), and of the Logscale plot of the relative cardinality
  of $\#\Lambda_{a, K}{(\Omega_0)}$ with respect to the set of all the integer
  vectors $\nu$ with $\norm{\nu}=K$. }
\label{fig3}
\end{figure}

Let us first consider any $n\geq 3$, $\Omega_0\ne (0,\ldots ,0)$, $1\le N+1\le M\in \mathbb N$, and let
\begin{eqnarray}\label{Lambdaa}
&&\Lambda_a{(\Omega_0)}:=\Big\{\nu\in{\mathbb Z}^n\setminus\{0\}:\ 
|\nu|\le M\,,\ 
|\Omega_0\cdot\nu|<a\Big\}\nonumber\\
&&  \Lambda_{a, K}{(\Omega_0)}:=\Big\{\nu\in \Lambda_a{(\Omega_0)}:\ |\nu|=K\Big\}\quad K=1\,,\ldots, M\nonumber\\
&& f_{N, M}:=\sum_{h\ne k\atop{N+1\le |h|+|k|\le M}} c_{hk}w^h z^k\,.
\end{eqnarray}
\begin{lemma}\label{LEM: resonant term bound}
The norm $\|f_*\|_{2R}$ of the function
$f_*$ defined in \eqref{fff} can be split as
\begin{eqnarray}\label{f*sum}
\|f_*\|_{2R}=\sum_{K=1}^M\|f_*^K\|_{2R}
\end{eqnarray}
where
\begin{eqnarray}\label{resonant term bound}
\|f^K_*\|_{2R}\le \#\Lambda_{a, K}{(\Omega_0)}\left(\frac{2R}{R_0}\right)^{\max\{K\,,\ N+1\}}\|f_{N, M}\|_{R_0}\,.
\end{eqnarray}
\end{lemma}
{In Figure \ref{fig3} (left panel) we represent the  cardinality $\#\Lambda_{a, K}{(\Omega_0)}$ for the vector $\Omega_0$ defined by the FPU case discussed in
Section 4 (see Eq. (\ref{Omega0FPU})) and the decrement,
for increasing $K$, of the ratio
between the cardinality of the set of all the integer
  vectors $\nu$ with $\norm{\nu}=K$ and $\#\Lambda_{a, K}{(\Omega_0)}$. 
  The cardinality is computed using Proposition \ref{cardinality3} presented
  in Subsection \ref{computationLambdan3}.}

\subsection{Proof of Lemma \ref{LEM: resonant term bound}}

\paragraph{\it Proof of Lemma \ref{LEM: resonant term bound}.} We write $f_*$ using the indices $h\in \mathbb N^n$ and $\nu:=h-k \in\mathbb Z^n\setminus\{0\}$:
\begin{eqnarray}\label{NEWf*}
f_*:=\sum_{\nu\in \Lambda_a{(\Omega_0)}}\sum_{h\in{\cal H}_*(\nu)}c_{hk}w^hz^{h-\nu}
\end{eqnarray}
where
  ${\cal H}_*(\nu)$ is defined by conditions
\begin{eqnarray}\label{H*(nu)}
{\cal H}_*(\nu):\quad \left\{
\begin{array}{ll}
h_j\ge 0\\
h_j\ge \nu_j\\
N+1\le|h|+|h-\nu|=\sum_{j=1}^n (2h_j-\nu_j)\le M
\end{array}
\right.
\end{eqnarray}
while,
using
\begin{eqnarray}\label{nu}
1\le |\nu|\le |h|+|k|\le M\,,
\end{eqnarray}
$\Lambda_a{(\Omega_0)}$ is as in \eqref{Lambdaa}.
We further decompose
\begin{eqnarray}
\Lambda_a{(\Omega_0)}=\bigcup_{K=1}^M\Lambda_{a, K}{(\Omega_0)}
\end{eqnarray}
and rewrite \eqref{NEWf*} as
\begin{eqnarray}\label{NEWNEWf*}
f_*:=\sum_{K=1}^M\sum_{\nu\in\Lambda_{a, K}{(\Omega_0)}}\sum_{h\in{\cal H}_*(\nu)}c_{hk}w^hz^{h-\nu}=\sum_{K=1}^Mf_*^K
\end{eqnarray}
with
\begin{eqnarray}
f_*^K:=\sum_{\nu\in\Lambda_{a, K}{(\Omega_0)}}\sum_{h\in{\cal H}_*(\nu)}c_{hk}w^hz^{h-\nu}\,.
\end{eqnarray}
From \eqref{NEWNEWf*}, it is immediate to recognize that
\begin{eqnarray}
\|f_*\|_{2R}=\sum_{K=1}^M\|f^K_*\|_{2R}\,.\end{eqnarray}
So, we are led to evaluate the norm of each $f^K_*$:
\begin{eqnarray}\label{norm f*KNEW}
\|f^K_*\|_{2R}=\sum_{\nu\in\Lambda_{a, K}{(\Omega_0)}}\sum_{h\in{\cal H}_*(\nu)}|c_{hk}|(2R)^{|h|+|h-\nu|}\,.
\end{eqnarray}
Using the definition of ${\cal H}_*(\nu)$ in \eqref{H*(nu)} and the inequality in \eqref{nu} with $k=h-\nu$, we have that the exponent of $2R$ in \eqref{norm f*KNEW} satisfies
\begin{eqnarray}
\max\{|\nu|\,,\ N+1\}\le |h|+|h-\nu|\le M\qquad \forall\ h\in {\cal H}_*(\nu)\,.
\end{eqnarray}
Hence, 
\begin{eqnarray}
\max\{K\,,\ N+1\}\le |h|+|h-\nu|\le M\qquad \forall\ h\in {\cal H}_*(\nu)\,,\ \forall\ \nu\in\Lambda_{a, K}{(\Omega_0)}\,.
\end{eqnarray}
This gives
\begin{eqnarray}
\|f^K_*\|_{2R}&\le& \left(\frac{2R}{R_0}\right)^{\max\{K\,,\ N+1\}}\sum_{\nu\in\Lambda_{a, K}{(\Omega_0)}}\sum_{h\in{\cal H}_*(\nu)}|c_{hk}|(2R_0)^{|h|+|h-\nu|}\nonumber\\
&\le& \#\Lambda_{a, K}{(\Omega_0)}\left(\frac{2R}{R_0}\right)^{\max\{K\,,\ N+1\}}\|f_{N, M}\|_{R_0}
\end{eqnarray}
having $\nu$--uniformly bounded the last summand from above by  $\|f_{N, M}\|_{R_0}$. $\quad \square$

\subsection{Computation of $\#\Lambda_{a, K}$ in the case $n=3$}\label{computationLambdan3}

Let $0<a\in \mathbb R$, 
 $0\ne K\in \natural$, $ \Omega_0:=(\Omega_{01}, \Omega_{02}, \Omega_{03})\in \mathbb R^3\setminus\{(0, 0, 0)\}$ be fixed.  Without loss of generality\footnote{
Up to renaming coordinates, one can always assume $|\Omega_{01}|\ge |\Omega_{02}|\ge |\Omega_{03}|$. As $(\Omega_{01}, \Omega_{02}, \Omega_{03})\ne (0, 0, 0)$, it must be
$\Omega_{01}\ne0$.
  As the set
 $\Lambda_{a,  K}$ in \eqref{Lambdaa} does not change 
changing $\Omega_0$ with $ \sigma\Omega_0$, where $\sigma=\pm 1$, taking $\sigma=\sign(\Omega_{01})$, one has \eqref{initial assumption}.},  assume that
\begin{eqnarray}\label{initial assumption}
 \Omega_{01}>0\,,\qquad  \Omega_{01}-|\Omega_{02}|\ge 0
\end{eqnarray}
For $(\sigma, \tau)\in \{(+, +), (+, -)\}$, we denote as
\begin{eqnarray}\label{asigmatauJKOLD}
a^{\sigma, \tau}_{J, K}(\Omega_0)&:=&\max\{
0\,,\ a^{\sigma, \tau}_{J, K,  1}(\Omega_0)\,,\ a^{\sigma, \tau}_{J, K,  2}(\Omega_0)
\}
\end{eqnarray}
where
\begin{eqnarray}\label{asigmatauJK}
a^{\sigma, \tau}_{J, K, 1}(\Omega_0)&:=&\left\{\begin{array}{lll}
\displaystyle
{\Omega_0\cdot(0, J, K-J)}
\quad &{\rm if}\ (\sigma, \tau)=(+, +)\\\\
\displaystyle
{\Omega_0\cdot (1, 1-J, K-J)}&{\rm if}\ (\sigma, \tau)=(+, -)
\end{array}
\right.\nonumber\\\nonumber\\
a^{\sigma, \tau}_{J, K, 2}(\Omega_0)&:=&\left\{\begin{array}{lll}
\displaystyle
{\Omega_0\cdot(1-J, -1, J-K)}\quad &{\rm if}\ (\sigma, \tau)=(+, +)\\\\
\displaystyle{\Omega_0\cdot (-J, 0, J-K)}&{\rm if}\ (\sigma, \tau)=(+, -)
\end{array}
\right.
\end{eqnarray}
{Put}
\begin{eqnarray}
{
w_{\sigma, \tau}}&{:=}&\left\{
\begin{array}{lll}
\displaystyle{+1}\quad &{{\rm if}\ (\sigma, \tau)=(+,+)}\\\\
\displaystyle{-1}  &{{\rm if}\ (\sigma, \tau)=(+,-)}
\end{array}
\right.
\end{eqnarray}
If 
 $\Omega_{01}-{w_{\sigma, \tau}}\Omega_{02}> 0$ {and $J=1$, $\ldots$, $K$}, put
 \begin{eqnarray}\label{nsigmatauaJK}
n^{\sigma, \tau, 0}_{a, J, K}(\Omega_0)&=&J\,,\qquad\qquad\qquad\quad n^{\sigma, \tau, 1}_{a, J, K}(\Omega_0)=
\left\{
\begin{array}{lll}
\left[\frac{a-a^{\sigma, \tau}_{J, K,  1}(\Omega_0)}{\Omega_{01}-{w_{\sigma\tau}}\Omega_{02}}\right]+1\quad &{\rm if}\  \frac{a-a^{\sigma, \tau}_{J, K,  1}(\Omega_0)}{\Omega_{01}-{w_{\sigma\tau}}\Omega_{02}}\notin {\mathbb Z}\\\\
\left[\frac{a-a^{\sigma, \tau}_{J, K,  1}(\Omega_0)}{\Omega_{01}-{w_{\sigma\tau}}\Omega_{02}}\right]&{\rm otherwise}
\end{array}
\right.
\nonumber\\
n^{\sigma, \tau, 2}_{a, J, K}(\Omega_0)&=&\left\{
\begin{array}{lll}
\left[\frac{a-a^{\sigma, \tau}_{J, K,  2}(\Omega_0)}{\Omega_{01}-{w_{\sigma\tau}}\Omega_{02}}\right]+1\quad &{\rm if}\ \frac{a-a^{\sigma, \tau}_{J, K,  2}(\Omega_0)}{\Omega_{01}-{w_{\sigma\tau}}\Omega_{02}}\notin {\mathbb Z}\\\\
\left[\frac{a-a^{\sigma, \tau}_{J, K,  2}(\Omega_0)}{\Omega_{01}-{w_{\sigma\tau}}\Omega_{02}}\right]&{\rm otherwise}
\end{array}
\right.\,,\nonumber\\
n^{\sigma, \tau, 3}_{a, J, K}(\Omega_0)&=&\left\{
\begin{array}{lll}
\left[\frac{a-a^{\sigma, \tau}_{J, K,  1}(\Omega_0)}{\Omega_{01}-{w_{\sigma\tau}}\Omega_{02}}\right]-\left[\frac{-a-a^{\sigma, \tau}_{J, K,  1}(\Omega_0)}{\Omega_{01}-{w_{\sigma\tau}}\Omega_{02}}\right]\ &\ \,  {\rm if}\quad 
\frac{a-a^{\sigma, \tau}_{J, K,  1}(\Omega_0)}{\Omega_{01}-{w_{\sigma\tau}}\Omega_{02}}\notin\mathbb Z\\\\
\max\left\{\left[\frac{a-a^{\sigma, \tau}_{J, K,  1}(\Omega_0)}{\Omega_{01}-{w_{\sigma\tau}}\Omega_{02}}\right]-\left[\frac{-a-a^{\sigma, \tau}_{J, K,  1}(\Omega_0)}{\Omega_{01}-{w_{\sigma\tau}}\Omega_{02}}\right]-1\,,0\right\}\ &\ \,  {\rm if}\quad 
\frac{a-a^{\sigma, \tau}_{J, K,  1}(\Omega_0)}{\Omega_{01}-{w_{\sigma\tau}}\Omega_{02}}\in\mathbb Z
\end{array}
\right.
\end{eqnarray}
Finally, {for $J=1$, $\ldots$, $K$}, put
 \begin{eqnarray}\label{msigmatauaJK}
  m^{\sigma, \tau}_{a, J, K}(\Omega_0)=\left\{
  \begin{array}{lll}
  \displaystyle0\quad &{\rm (i)}\ {\rm if}\ a_{J, K}^{\sigma, \tau}(\Omega_0)>0\ {\rm and}\ a\le a_{J, K}^{\sigma, \tau}\\\\
  \min\big\{n^{\sigma, \tau, s}_{a, J, K}(\Omega_0)\big\}_{s=0\,,1\,,2\,,3}&{\rm (ii)}\ {\rm if}\ a> a_{J, K}^{\sigma, \tau}(\Omega_0)\ {\rm and}\ \Omega_{01}-{w_{\sigma\tau}}\Omega_{02}>0\\\\
    \displaystyle\ J&{\rm (iii)}\ {\rm if}\ a> a_{J, K}^{\sigma, \tau}(\Omega_0)\ {\rm and}\ \Omega_{01}-{w_{\sigma\tau}}\Omega_{02}=0
  \end{array}
  \right.
  \end{eqnarray}
  {and
  \begin{eqnarray}\label{m}
  m_{a, J, K}(\Omega_0):=\left\{
  \begin{array}{lll}
 \displaystyle 1\quad &{\rm if}\ J=0\ {\rm and}\ |\Omega_{03}|K<a\\\\
  \displaystyle 0\quad &{\rm if}\ J=0\ {\rm and}\ |\Omega_{03}|K\ge a\\\\
  \displaystyle \sum_{(\sigma, \tau)\in \{(+,+), (+, -)\}}m^{\sigma,\tau}_{a, J, K}(\Omega_0)+m^{\sigma, \tau}_{a, J, K}(\Omega^{\rm rev}_0)\quad &{\rm if}\ J=1, \ldots, K\\
  \end{array}
  \right.
  \end{eqnarray}}
  {where $ \Omega^{\rm rev}_0:=(\Omega_{01}, \Omega_{02}, -\Omega_{03})$.}
\begin{proposition}\label{cardinality3} Let $n=3$, $a>0$, 
  $1 \leq K\in \natural$, $ \Omega_0\in \mathbb R^3$ satisfying \eqref{initial assumption}. Then\footnote{The symbol  $\delta_{m, n}$ is the Kronecker symbol.}
  \begin{eqnarray}\label{cardinality 31}{\#\Lambda_{a, K}(\Omega_0)=\sum_{J=0}^K(2-\delta_{J, K})m_{a, J, K}(\Omega_0)}\,.
\end{eqnarray}
\end{proposition}

\subsection{Proof of Proposition \ref{cardinality3}}
To prove  Proposition \ref{cardinality3} we decompose
\begin{eqnarray}\label{decomposition1}
\Lambda_{a, K}(\Omega_0)=\bigcup_{J=0}^K\Lambda_{a, J, K}(\Omega_0)\end{eqnarray}
and, further
\begin{eqnarray}\label{decomposition2}
\Lambda_{a, J, K}(\Omega_0)=\Lambda^{(-)}_{a, J, K}(\Omega_0)\cup\Lambda^{(0)}_{a,  K}(\Omega_0)\cup \Lambda^{(+)}_{a, J, K}(\Omega_0)\end{eqnarray}
and, finally, \begin{eqnarray}\label{decomposition3}
\Lambda^{(+)}_{a, J, K}(\Omega_0)=\bigcup_{(\sigma, \tau)\in\{+, -\}^2}\Lambda^{\sigma, \tau}_{a, J, K}(\Omega_0)\quad {\forall\ J=1\,,\ldots, K}\end{eqnarray}
where, for any
$0\le J\le K\in \mathbb N$, 
\begin{eqnarray}\label{set Lambda}
\Lambda_{a, J, K}(\Omega_0)&:=&\Big\{\nu=(\nu_1, \nu_2, \nu_3)\in\Lambda_{a, K}(\Omega_0):\ |\nu_1|+|\nu_2|=J\Big\}\nonumber\\
&=&\Big\{\nu=(\nu_1, \nu_2, \nu_3)\in {\mathbb Z}^3\setminus\{0\}:\ \left|
\Omega_0\cdot \nu
\right|<a\,,\ |\nu_1|+|\nu_2|=J\,,\nonumber\\
&&|\nu_1|+|\nu_2|+|\nu_3|=K\Big\}\,.\end{eqnarray}
and
\begin{eqnarray}\label{Lambda+0-NEW}
&&\Lambda^{(-)}_{a, J, K}(\Omega_0):=\Lambda_{a, J, K}(\Omega_0)\cap\{\nu_3<0\}\,,\ \Lambda^{(0)}_{a,  K}(\Omega_0):=\Lambda_{a, K, K}(\Omega_0)\cap\{\nu_3=0\}\nonumber\\
&&\Lambda^{(+)}_{a, J, K}(\Omega_0):=\Lambda_{a, J, K}(\Omega_0)\cap\{\nu_3>0\}
\end{eqnarray}
and, {for $J=1$, $\ldots$, $K$},
\begin{eqnarray}\label{Lambdasigmatau}
\Lambda^{+, +}_{a, J, K}(\Omega_0)&:=&\Lambda^{(+)}_{a, J, K}(\Omega_0)\cap\Big\{ \nu_1\ge 0\,,\ \nu_2>0\Big\}\nonumber\\
\Lambda^{+, -}_{a, J, K}(\Omega_0)&:=&\Lambda^{(+)}_{a, J, K}(\Omega_0)\cap\Big\{\nu_1>0\,,\ \nu_2\le 0\Big\}\nonumber\\
\Lambda^{ -, +}_{a, J, K}(\Omega_0)&:=&\Lambda^{(+)}_{a, J, K}(\Omega_0)\cap\Big\{\nu_1<0\,,\ \nu_2\ge 0\Big\}\nonumber\\ \Lambda^{-, -}_{a, J, K}(\Omega_0)&:=&\Lambda^{(+)}_{a, J, K}(\Omega_0)\cap\Big\{ \nu_1\le 0\,,\ \nu_2<0\Big\}
\end{eqnarray}

\begin{lemma}\label{symmetries}\item[{\rm(i)}]
 $\Lambda^{(-)}_{a, J, K}(\Omega_0)$ and $\Lambda^{(+)}_{a, J, K}(\Omega_0)$ have the same cardinality;
\item[{\rm(ii)}] the cardinality of $\Lambda^{(0)}_{a, K}(\Omega_0)$ can be deduced from the one of $\Lambda^{(-)}_{a, J, K}(\Omega_0)$ or $\Lambda^{(+)}_{a, J, K}(\Omega_0)$ taking, formally\footnote{According to the definitions in \eqref{Lambda+0-NEW},  $\Lambda^{(-)}_{a, J, K}(\Omega_0)$ and $\Lambda^{(+)}_{a, J, K}(\Omega_0)$ 
are meaningful only when $J<K$.}, $J=K$;
\item[{\rm(iii)}]  $\#\Lambda^{-, -}_{a, J, K}(\Omega_0)=\#\Lambda^{+, +}_{a, J, K}(\Omega^{\rm rev}_0)$ and
$\#\Lambda^{-, +}_{a, J, K}(\Omega_0)=\#\Lambda^{+, -}_{a, J, K}(\Omega^{\rm rev}_0)$.
\end{lemma}

\paragraph{\it Proof of Lemma  \ref{symmetries}.} (i) The bijection
$
(\nu_1, \nu_2, \nu_3)\to (-\nu_1, -\nu_2, -\nu_3)
$
carries $\Lambda^{(-)}_{a, J, K}$ to $\Lambda^{(+)}_{a, J, K}$. 
 Item (ii) follows from the rewrites
\begin{eqnarray}
\Lambda^{(-)}_{a, J, K}&=&\Big\{(\nu_1, \nu_2)\in {\mathbb Z}^2:\  -a<\Omega_{01}\nu_1+\Omega_{02}\nu_2-(K-J)\Omega_{03}<a\,,\ |\nu_1|+|\nu_2|=J\Big\}\nonumber\\
&&\times\{-(K-J)\}\nonumber\\
\Lambda^{(0)}_{a, K\phantom{, J}}&=&\Big\{(\nu_1, \nu_2)\in {\mathbb Z}^2:\  -a<\Omega_{01}\nu_1+\Omega_{02}\nu_2<a\,,\ |\nu_1|+|\nu_2|=K\Big\}\times\{0\}\nonumber\\
\Lambda^{(+)}_{a, J, K}&=&\Big\{(\nu_1, \nu_2)\in {\mathbb Z}^2:\  -a<\Omega_{01}\nu_1+\Omega_{02}\nu_2+(K-J)\Omega_{03}<a\,,\ |\nu_1|+|\nu_2|=J\Big\}\nonumber\\
&&\times\{K-J\}\nonumber\\
\end{eqnarray}
(iii) Combining the bijection $
(\nu_1, \nu_2)\to (-\nu_1, -\nu_2)
$
with the change $\Omega_{03}\to -\Omega_{03}$
carries
 $\Lambda^{-, +}_{a, J, K}(\Omega_0)$,   to
 $\Lambda^{+, -}_{a, J, K}(\Omega_0^{\rm rev})$,   and $\Lambda^{-, -}_{a, J, K}(\Omega_0)$ to $\Lambda^{+,+}_{a, J, K}(\Omega^{\rm rev}_0)$.
$\quad \square$

 \begin{lemma}\label{cardinality32}
 Let $0<a\in \mathbb R$, {$K\ne 0$,}
$0\le J\le K\in \mathbb N$, $ \Omega_0:=(\Omega_{01}, \Omega_{02}, \Omega_{03})\in \mathbb R^3\setminus\{(0, 0, 0)\}$,
  $(\sigma, \tau)\in \{(+, +), (+, -)\}$. Assume \eqref{initial assumption}. Then:{
  \begin{itemize}
  \item[$a)$] the cardinality of $\Lambda^{(+)}_{a, 0, K}$ is the number $m_{a, 0, K}$ in \eqref{m};
  \item[$b)$] the cardinality of $\Lambda^{\sigma, \tau}_{a, J, K}$ is the number $m^{\sigma, \tau}_{a, J, K}$ in \eqref{msigmatauaJK}.
  \end{itemize}
  }
\end{lemma}
\begin{remark}\rm By definition \eqref{asigmatauJKOLD}, \eqref{asigmatauJK}, $a^{\sigma, \tau}_{J, K}$ is greater of equal than $a_{J, K, 1}^{\sigma, \tau}$ and $a_{J, K, 2}^{\sigma, \tau}$ both. Thus,
if $a> a_{J, K}^{\sigma, \tau}$, the $n^{\sigma, \tau, s}_{a, J, K}$ with $s\in\{0, 1, 2\}$ are all positive. However, $\Lambda^{\sigma, \tau}_{a, J, K}$ may still be empty, if it happens $n^{\sigma, \tau, 3}_{a, J, K}(\Omega_0)=0$.
\end{remark}

\paragraph{\it Proof of Lemma \ref{cardinality32}.}{$a)$ is obvious. We discuss $b)$.} 
The cardinality of the sets $\Lambda^{+, +}_{a, J, K}$, $\Lambda^{+, -}_{a, J, K}$ coincides with the one of the index sets
\begin{eqnarray}\label{JsigmatauaJK}
{\rm J}^{+, +}_{a, J, K}&:=&\Big\{
i\in \mathbb Z:\ 0\le i\le J-1\,,\ -a<\Omega_{01}i+\Omega_{02}(J-i)+\Omega_{03}(K-J)<a\Big\}\nonumber\\
{\rm J}^{+, -}_{a, J, K}&:=&\Big\{
i\in \mathbb Z:\ 
 0\le i\le J-1\,,\ -a<\Omega_{01}(i+1)-\Omega_{02}(J-i-1)+\Omega_{03}(K-J)<a\Big\}\nonumber\\
\end{eqnarray}
Indeed, one has
\begin{eqnarray}
\Lambda^{+, +}_{a, J, K}&:=&\Big\{
(i\,,\  J-i\,,\ K-J):\ i\in {\rm J}^{+, +}_{a, J, K}\Big\}\nonumber\\
\Lambda^{+, -}_{a, J, K}&:=&\Big\{
(i+1\,,\  -(J-i-1)\,,\  K-J):\ i\in {\rm J}^{+, -}_{a, J, K}\Big\}\,.
\end{eqnarray}
We then compute $\#{\rm J}^{+, +}_{a, J, K}$, $\#{\rm J}^{+, -}_{a, J, K}$. 
We consider the case $\Omega_{01}-{w_{\sigma\tau}}\Omega_{02}>0$.
In such case,
${\rm J}^{\sigma, \tau}_{a, J, K}$ is the set of integer numbers in the interval \begin{eqnarray}{\cal J}^{\sigma, \tau}_{a, J, K}:=\bigcap_{s=0}^4{\cal J}^{\sigma, \tau, s}_{a, J, K}\end{eqnarray} where
\begin{eqnarray}\label{cal J intervals}
&&{\cal J}^{\sigma, \tau, 0}_{a, J, K}:=\big[0\,,J-1\big]\,,\qquad {\cal J}^{\sigma, \tau, 1}_{a, J, K}:=\left[0\,,\ \frac{a-\gamma^{\sigma, \tau}_{J, K}}{\Omega_{01}-{w_{\sigma\tau}}\Omega_{02}}\right)\nonumber\\
&&{\cal J}^{\sigma, \tau, 2}_{a, J, K}:=\left(\frac{-a-\gamma^{\sigma, \tau}_{J, K}}{\Omega_{01}-{w_{\sigma\tau}}\Omega_{02}}\,,\ J-1\right]\,,\quad {\cal J}^{\sigma, \tau, 3}_{a, J, K}:=\left(\frac{-a-\gamma^{\sigma, \tau}_{J, K}}{\Omega_{01}-{w_{\sigma\tau}}\Omega_{02}}\,,\ \frac{a-\gamma^{\sigma, \tau}_{J, K}}{\Omega_{01}-{w_{\sigma\tau}}\Omega_{02}}\right)
\end{eqnarray}
with
\begin{eqnarray}
\gamma^{+, +}_{J, K}:=\Omega_{02}J+
\Omega_{03}(K-J)\,,\quad 
\gamma^{+, -}_{J, K}:=\Omega_{01}-\Omega_{02}(J-1)+\Omega_{03}(K-J)\,.\end{eqnarray}
The sets
${\cal J}^{\sigma, \tau, s}_{a, J, K}$ are non--empty, as  intervals in $\mathbb R$, if and only if
\begin{eqnarray}\label{existence conditions NEW}
a>\max\{a^{\sigma, \tau}_{J, K,  1}(\Omega_0)\,,\ a^{\sigma, \tau}_{J, K,  2}(\Omega_0)\}
\end{eqnarray}
where $a^{\sigma, \tau}_{J, K,  1}$, $a^{\sigma, \tau}_{J, K,  2}$ are as in \eqref{asigmatauJK}. This proves (i), in the case $\Omega_{01}-w_{\sigma\tau}\Omega_{02}>0$.
Moreover, by definition, ${\cal J}^{\sigma, \tau}_{a, J, K}$ is one of the of the ${\cal J}^{\sigma, \tau, s}_{a, J, K}$ in \eqref{cal J intervals} and is part of all of them. Thus, the amount of integer numbers in ${\cal J}^{\sigma, \tau}_{a, J, K}$ coincides with the minimum among the amounts of integer numbers in ${\cal J}^{\sigma, \tau, s}_{a, J, K}$, which are precisely the numbers $n^{\sigma, \tau, s}_{a, J, K}$ in \eqref{nsigmatauaJK}. We have thus proved (ii). The proof of item (iii) and the completion of (i) when $\Omega_{01}-{w_{\sigma\tau}}\Omega_{02}=0$  are particularizations  of the formulae \eqref{JsigmatauaJK}. $\qquad \square$

\paragraph{\it Proof of Proposition \ref{cardinality3}.} Using the decompositions \eqref{decomposition1}, \eqref{decomposition2}, \eqref{decomposition3}, Lemma \ref{symmetries}, Lemma \ref{cardinality32}
\begin{eqnarray}\#\Lambda_{a, K}{(\Omega_0)}&=&\sum_{J=0}^K\#\Lambda_{a, J, K}{(\Omega_0)}\nonumber\\
&=&
\sum_{J=0}^K\Big(\#\Lambda^{(-)}_{a, J, K}{(\Omega_0)}+\#\Lambda^{(0)}_{a, J, K}{(\Omega_0)}+\#\Lambda^{(+)}_{a, J, K}{(\Omega_0)}\Big)\nonumber\\
&=&\sum_{J=0}^{K}
(2-\delta_{J, K})\#\Lambda^{(+)}_{a, J, K}{(\Omega_0)}\nonumber\\
&=&
{2\#\Lambda^{(+)}_{a, 0, K}{(\Omega_0)}+}
\sum_{J={1}}^{K}
\sum_{(\sigma, \tau)\in\{(+, +), (+, -)\}}
(2-\delta_{J, K})\Big(\#\Lambda^{\sigma, \tau}_{a, J, K}(\Omega_0)+\#\Lambda^{\sigma, \tau}_{a, J, K}(\Omega^{\rm rev}_0)\Big)
\nonumber\\
&=&{2m_{a, 0, K}(\Omega_0)+}\sum_{J={1}}^{K}
\sum_{(\sigma, \tau)\in\{(+, +), (+, -)\}}
(2-\delta_{J, K})\Big(
m^{\sigma, \tau}_{a, J, K}(\Omega_0)+m^{\sigma, \tau}_{a, J, K}(\Omega^{\rm rev}_0)\Big)\nonumber\\
&=&{\sum_{J=0}^{K}
(2-\delta_{J, K})m_{a, J, K}(\Omega_0)}\,.
\end{eqnarray}
{In the last step, we have used the definitions \eqref{m} and that $K\ne 0$, whence $\delta_{0, K}=0$. $\quad \square$}


 \section{Appendix: non-stationary phase approach for non--resonant
 normal forms of Nekhoroshev's theory}\label{argumentforAA}

 Let us explain the idea that we exploit in this paper in 
the simpler framework of the action--angle analytic Hamiltonians having the form
\begin{equation}\label{action angle}
    H_0(I, \varphi)=h(I)+\varepsilon f(I, \varphi)
\end{equation}
where $I\in V\subseteq {\mathbb R}^n$, $V$ open and connected, $\varphi\in {\mathbb T}^n$,  $\varepsilon$ is a small positive number, and $h$ satisfies a generic condition known as ``steepness''. Within the proof of Nekhoroshev's theorem, the Hamiltonian (\ref{action angle}) is conjugate by a close--to--the--identity canonical transformation
to a non-resonant normal form:
\begin{equation}\label{normalformAA}
  H(I,\varphi)=h(I)+\varepsilon g(I)+\varepsilon r(I,\varphi)
\end{equation}  
defined for all $I\subseteq D \subset V$ where
$D\subseteq {\Bbb R}^n$ is a suitable open set, where the frequency vector
$\omega_\varepsilon(I)=\nabla_I (h(I)+\varepsilon g(I))$ satisfies the non-resonance condition:
\begin{equation}
  \norm{k\cdot \omega(I)}\geq \gamma_\varepsilon \, \ \  \forall\ k\in {\mathbb Z}^n\,: \quad 0<\sum_{j=1}^n|k_j|\le N\,,
\end{equation}
for suitable $\gamma_\varepsilon >0,N\sim 1/\varepsilon^c$ ($c>0$ depends
only on the number of the degrees of freedom and on the values of
the steepness indices of $h(I)$, see \cite{guzzoCB16});
the Fourier norm of the
remainder $r(I,\varphi)=\sum_{k\in {\Bbb Z}^n} r_k(I)e^{i k\cdot \varphi}$ becomes
exponentially small for increasing values of $N$.

Consider any solution of the Hamilton equation of (\ref{normalformAA})
with initial conditions with $I(0)\in D$. The classical way to bound the motion of the action variables for the normal form Hamiltonian (\ref{normalformAA}), which is also used within the proof of Nekhoroshev theorem, works as follows: for all the times $t$ such that  $I(t)\in D$, the following inequality is valid:
 \begin{equation}\label{a priori boundAA}
    |I_j(t)-I_j(0)|\le 
\sum_{k\in {\Bbb Z}^n}\norm{k_j} \left|\int_0^t r_k(I(\tau))e^{\ii k\cdot \varphi(\tau)}d\tau
\right|
    \leq  |t|\left (\sum_{k\in {\Bbb Z}^n} \norm{k_j}  \sup_{D}\norm{r_k(I)}\right )   .
    \end{equation}
We obtain an alternative estimate by individually revising the estimate of each term
 \begin{equation}\label{intNS}
   \left|\int_0^t r_k(I(\tau))e^{\ii k\cdot \varphi(\tau)}d\tau
   \right|
 \end{equation}
 depending on whether
\begin{equation}\label{smalldivaham}
  \inf_{I\in B(I_*,\rho)}\norm{k\cdot \omega_\varepsilon(I)} \geq a
\end{equation}
where $I_*=I(0)$ and $a,\rho>0$ are suitable parameters (they need to satisfy some standard consistency conditions), $B(I_*,\rho )\subseteq D$ is a real
ball of radius $\rho$ centered at $I_*$; notice that the parameter $a$ is not
necessarily set equal to $\gamma_\varepsilon$.

For the integer vectors $k$ such that inequality (\ref{smalldivaham}) is not
satisfied, we proceed as in
 the traditional way, and obtain:
\begin{equation}\label{ineqrktrad}
\left|\int_0^t r_k(I(\tau))e^{\ii k\cdot \varphi(\tau)}d\tau
  \right| \leq |t|\sup_{D}\norm{r_k(I)} .
\end{equation}
An improvement is instead obtained if inequality (\ref{smalldivaham}) is
satisfied, using a representation valid for integrals with non-stationary
 phases (see for example \cite{asymptoticintegrals}). Indeed,
 we have:
\begin{eqnarray}\label{NSestimateAA}
    \int_0^t r_{k}(I(\tau))
    e^{ \ii k \cdot \varphi(\tau)}d\tau    &=&
    \int_0^t {  r_{k}(I(\tau))\over \ii k \cdot \omega_\varepsilon (I(\tau)) }
    \Big [ \ii k\cdot \dot \varphi (\tau)- \ii \varepsilon k \cdot {\partial r
        \over \partial I}(I(\tau),\varphi(\tau)) \Big ]
 e^{ \ii k \cdot \varphi(\tau)}d\tau 
\cr
&=& 
 \int_0^t {  r_{k}(I(\tau))\over \ii k \cdot \omega_\varepsilon (I(\tau)) }{d\over d\tau}e^{ \ii k \cdot \varphi(\tau)}d\tau 
 - \int_0^t {  \varepsilon r_{k}(I(\tau))\over  k \cdot \omega_\varepsilon (I(\tau)) }  k \cdot {\partial r
        \over \partial I}(I(\tau),\varphi(\tau) )
 e^{ \ii k \cdot \varphi(\tau)}d\tau 
\cr
&=&    {r_{k}(I(\tau))\over  \ii k \cdot \omega_\varepsilon (I(\tau))}
 e^{ \ii k \cdot \varphi(\tau)}\Big \vert_0^t
 -  \int_0^t \left \{ {r_{k}(I)\over  \ii k \cdot \omega_\varepsilon(I)},
 \varepsilon r\right \}_{\vert (I(\tau),\varphi(\tau))}    e^{ \ii k \cdot \varphi(\tau)}d\tau\cr
 &&
 - \int_0^t {  \varepsilon r_{k}(I(\tau))\over  k \cdot \omega_\varepsilon (I(\tau)) }   k \cdot {\partial r
        \over \partial I}(I(\tau),\varphi(\tau) )
 e^{ \ii k \cdot \varphi(\tau)}d\tau
 \end{eqnarray}
From Eq. (\ref{NSestimateAA}) we get indeed an improvement of
inequality (\ref{ineqrktrad}) for suitably long time $\norm{t}$. In fact,
both integrals in the right hand side of Eq. (\ref{NSestimateAA}) have
integrands of second order in the remainder $r(I,\varphi)$, while
the first term, which is only linear in the remainder $r(I,\varphi)$, is
just a boundary term (its estimate is not multiplied by $\norm{t}$). 

The full exploitation of the representation (\ref{NSestimateAA}) is
obtained when all the Fourier harmonics $r_k(I)$ of the remainder
are explicitly known (or bounded within an interval) up to
an order $\sum_i \norm{k_i}  \leq M$, with $M>N$: for the
$k\in {\Bbb Z}$ with $N<\sum_i \norm{k_i}  \leq M$ we use
(\ref{intNS}) or the representation (\ref{NSestimateAA}) according
to the value of $\alpha_k=\inf_D \norm{k\cdot \omega(I)}$; 
for the $k\in {\Bbb Z}$ with $ \sum_i \norm{k_i}> M$ we use
(\ref{intNS}). 
\vskip 0.4 cm
\noindent
{\bf Acknowledgments.} The authors acknowledge A. Giorgilli, U. Locatelli and M. Sansottera for the availability of the algebraic manipulator {\grm Qr'onoc}.

\end{document}